\documentclass{svjour3}
\usepackage{amssymb}

\usepackage{epic}
\usepackage{amsmath}
\usepackage{amsfonts}
\usepackage{pstricks}

\usepackage{float}
\usepackage{subfig}
\usepackage{graphicx}
\usepackage{multirow}
\usepackage{color}
\usepackage[subfigure]{tocloft}
\usepackage{subfig}
\usepackage[utf8]{inputenc}
\usepackage{soul}
\newtheorem{propr}{Property}
\usepackage{hyperref}
\hypersetup{colorlinks=true,linkcolor=blue,urlcolor=blue}
\def\argmin{\mathop{\text{argmin}}}
\def\tw{\widetilde{\omega}}

\journalname{Journal of Scientific Computing}

\begin{document}
\title{High order weighted extrapolation for boundary conditions for
  finite difference methods on complex domains with Cartesian meshes}
\titlerunning{High order boundary extrapolation}
\author{A.~Baeza \and P.~Mulet \and D.~Zorío}
\date{}

\institute{
Departament de Matemàtica Aplicada,  Universitat de
  València    (Spain); emails: antonio.baeza@uv.es, mulet@uv.es, david.zorio@uv.es. This research was partially supported by Spanish MINECO grants
   MTM2011-22741 and MTM2014-54388-P.
}

\leavevmode\thispagestyle{empty}

\noindent This version of the article has been accepted for publication, after a peer-review process, and is subject to Springer Nature’s AM terms of use, but is not the Version of Record and does not reflect post-acceptance improvements, or any corrections. The Version of Record is available online at: \url{https://doi.org/10.1007/s10915-016-0188-7}

\newpage

\maketitle

\begin{abstract}
  The design of numerical boundary conditions is a
    challenging problem
  that has been tackled in different ways depending on the nature of
  the problem and the numerical scheme used to solve it. In this paper
  we present a new weighted extrapolation technique which entails an
  improvement with respect to the technique that was developed in
  \cite{BaezaMuletZorio2015}. This technique is based on the
  application of a variant of the Lagrange extrapolation through the
  computation of weights capable of detecting regions with
  discontinuities. We also present a combination of the above
  technique with a least squares approach in order to stabilize the
  scheme in some cases where Lagrange extrapolation can turn the
  scheme mildly unstable. We show that this combined extrapolation
  technique can tackle discontinuities more robustly than the
  procedure introduced in \cite{BaezaMuletZorio2015}.

\end{abstract}

\section{Introduction}
There are many physical phenomena that are described
  by hyperbolic
conservation laws. The lack of knowledge about the analytic solution
of most of them motivated the development of numerical methods to
approximate their solutions in the last decades.

In order to solve them accurately, the schemes evolved to improve the
accuracy and to capture well
the discontinuities in weak solutions. Due either to the nature of
the problems or to the impossibility of handling unbounded domains on
a computer numerical boundary conditions have to be imposed.

Despite high order numerical methods have been widely
  studied and
developed, usually the boundary treatment has consisted in low order
extrapolations, decreasing the global order of the scheme and the
quality of the simulation results.

In this  paper we present a weighted extrapolation technique for the
application of numerical boundary conditions in numerical schemes
using Cartesian meshes, which represents an improvement of our previous work 
\cite{BaezaMuletZorio2015}, where a Boolean
node selection with a thresholding criterion was used. 
In this paper we will design a fully weighted boundary extrapolation, with
scale independent and dimensionless weights. We will combine
these techniques with a least squares extrapolation to
overcome some stability issues that may appear in some
smooth problems with complex boundary shape. 

This work can be understood as an extension
to a higher order of accuracy  of \cite{Sjogreen},
where second order Lagrange  
extrapolation with slope limiters is used to extrapolate to a single 
ghost cell. A related approach for data extrapolation at the boundary 
with high order is proposed in \cite{TanShu} and \cite{TanWangShu}, 
where the authors develop an extrapolation technique for inflow boundaries
(known as inverse Lax-Wendroff) which is based on
using the equations to compute normal derivatives and approximate ghost 
values by means of a Taylor expansion that involves the 
Dirichlet boundary conditions and the computed derivatives,
whereas ghost cell values at outflow boundaries
are approximated by means of combinations of Lagrange and WENO 
extrapolation and least-squares approximation.
The aforementioned technique relies on weights that are not adimensional and scale dependent. 
Later on, the authors of \cite{FilbertYang} introduced a modification 
that overcame the dependency on the scale.

The organization of the paper is the following: In section \ref{seq}
we present the equations and the numerical methods that we
consider in this paper. The details of the procedure for meshing complex domains
with Cartesian meshes are explained in section \ref{sml},
where we
also expound how we perform extrapolations at the boundary. In section
\ref{sep} we present a proposal for weighted extrapolation
and introduce several procedures that can be used to design the weights. Some
numerical results that are
obtained with this methodology are presented in section \ref{srn},
with some simple tests in 1D to illustrate the correct behavior of the
proposed techniques and some more complex ones in 2D. Finally, some
conclusions are drawn in section \ref{scn}.

\section{Numerical scheme}

\label{seq}
We will consider along this paper hyperbolic systems of $m$
two-dimensional conservation laws
\begin{equation}\label{eq:hypsis}
u_t+f(u)_x+g(u)_y=0,\quad u=u(x,y,t),
\end{equation}
defined on an open and bounded spatial domain
$\Omega\subseteq\mathbb{R}^2$, with
Lipschitz boundary $\partial\Omega$, a finite union of piece-wise smooth
curves, $u:\Omega\times\mathbb{R}^+\rightarrow
\mathbb{R}^m$, and fluxes $f,g:\mathbb{R}^m\rightarrow\mathbb{R}^m$.
  These equations are supplemented with an initial condition, $u(x,
  y, 0)=u_0(x, y)$, $u_0:\Omega\rightarrow\mathbb{R}^m$,
  and different boundary conditions that may vary depending on the
  problem.

  The techniques that will be expounded in this paper can be used in a
  more general family of numerical schemes; however,
  we use here Shu-Osher's finite difference conservative methods
  \cite{ShuOsher1989} with a
  WENO5 (\textit{Weighted Essentially
    Non-Oscillatory}) \cite{JiangShu96} spatial reconstruction,
  Donat-Marquina's
  flux-splitting \cite{DonatMarquina96}  and the RK3-TVD ODE
  solver \cite{ShuOsher89}. This
  combination of techniques was proposed in \cite{MarquinaMulet03}.

\section{Meshing procedure}
\label{sml}
The mesh is defined through the following items: a reference vertical
line,
  $x=\overline{x}$, a horizontal line $y=\overline{y}$ and two
  positive values $h_x>0$
  and  $h_y>0$, the horizontal and vertical spacings of the mesh, so
  that the vertical lines in the mesh are determined by:
$x=x_r:=\overline{x}+rh_x$, $r\in\mathbb{Z}$ and
the horizontal ones by $ y=y_s:=\overline{y}+s h_y$,
$s\in\mathbb{Z}$.
The computational domain is given by the intersections of the horizontal and vertical mesh lines
that belong to the physical domain:
$$\mathcal{D}:=\left\{(x_r, y_s):\hspace{0.2cm}
(x_r,y_s)\in\Omega,\hspace{0.3cm}r,s\in\mathbb{Z}\right\}=
(\overline{x}+h_x\mathbb{Z})\times
(\overline{y}+h_y\mathbb{Z})\cap\Omega.$$

Each point $(x_r, y_s)$ determines a rectangular cell by:
\begin{equation*}
  [x_{r}-\frac{h_x}{2}, x_{r}+\frac{h_x}{2}]\times
  [y_{s}-\frac{h_y}{2}, y_{s}+\frac{h_y}{2}].
\end{equation*}

 In order to advance in time using WENO schemes of order $2k-1$, 
$k$ additional cells are needed at both sides of
each horizontal and vertical mesh line. If these additional
 cells fall outside the domain they are usually named \textit{ghost cells} and, in
terms of their centers, are given by:
$$\mathcal{GC}:=\mathcal{GC}_x\cup\mathcal{GC}_y,$$
where
$$\mathcal{GC}_x:=\left\{(x_r, y_s):\hspace{0.2cm}0<d\left(x_r,\hspace{0.1cm}\Pi_x\left(\mathcal{D}\cap\left(\mathbb{R}\times\{y_s\}\right)\right)\right)\leq kh_x,\hspace{0.3cm}r,s\in\mathbb{Z}\right\},$$
$$\mathcal{GC}_y:=\left\{(x_r,
  y_s):\hspace{0.2cm}0<d\left(y_s,\hspace{0.1cm}\Pi_y\left(\mathcal{D}\cap\left(\{x_r\}\times\mathbb{R}\right)\right)\right)\leq
  kh_y,\hspace{0.3cm}r,s\in\mathbb{Z}\right\},$$
with $\Pi_x$ and $\Pi_y$ being the projections on the respective
coordinates and,
$$d(a, B):=\inf\{|b-a|:\hspace{0.1cm}b\in B\},$$
for given   $a\in\mathbb{R}$ and
$B\subseteq\mathbb{R}$.
Notice that
$d(a,\emptyset)=+\infty$, since, by convention,
$\inf\emptyset=+\infty$.

There are many ways to compute the ghost cell values that are required for 
the numerical scheme to advance in time. In order 
to maintain the accuracy and stability of the numerical method
many issues are to be taken into account.

First of all and regarding the node organization and selection it seems reasonable to use
the lines normal to the boundary as directions for extrapolation, since many physical boundary conditions (e.g. Dirichlet) 
use the normal lines in their definition.
It seems thus reasonable that the
extrapolation at a certain ghost cell $P$ be based on the
physical conditions at its nearest boundary point, $P_0$, a solution of
 $$\|P-P_0\|_2=\min\{\|P-B\|_2:\quad B\in\partial\Omega\}.$$
 Uniqueness of $P_0$
holds whenever $P$ is close enough to the boundary, so
 we will henceforth denote $N(P)=P_0$.
 The line determined by $P$ and $N(P)$ is the  line normal to
  $\partial\Omega$ at  $N(P)$ if $\partial\Omega$ is
 differentiable at $P_0$.
 In a general setting nodes inside the domain do not fall along boundary     
normals and therefore a prior interpolation process is to be made at some nodes lying at the normal. For simplicity we will not distinguish between
interpolation or extrapolation, and the second term will be routinely used for both procedures. The motivation of the choice of the nodal disposition that we will  introduce next and further details about it can be found in \cite{BaezaMuletZorio2015}.

If a numerical scheme of order $r+1$ is used to numerically solve the 
equations and we wish to design a procedure for numerical boundary conditions 
that formally preserve the accuracy in the resulting
scheme it is reasonable to use extrapolation of the same order as
the numerical method (or higher) at the ghost cells.

\begin{figure}[htb]
      \centering
      \begin{tabular}{cc}
      \includegraphics[width=0.47\textwidth]{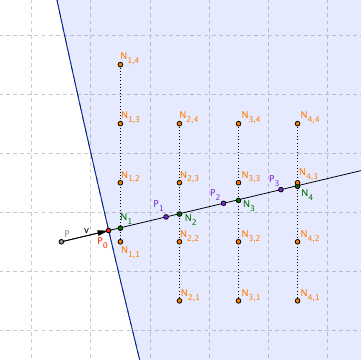} &
      \includegraphics[width=0.47\textwidth]{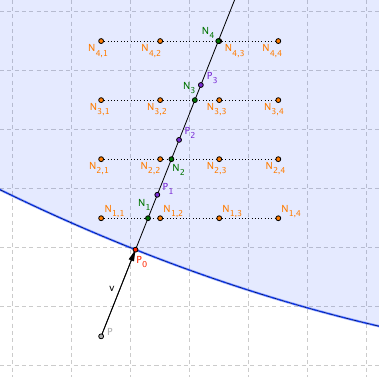}\\
      (a) & (b)
    \end{tabular}
    \caption{Examples of choice of stencil. We use the stencil
      $\mathcal{S}(P)=\{N_1, N_2, N_3, N_4\}$ in case of outflow
      boundary and conditions and the stencil $\mathcal{S}(P)=\{P_0,
      P_1, P_2, P_3\}$ in case of Dirichlet boundary conditions.}
      \label{fig:CxyW}
    \end{figure}

We proceed in a similar fashion as in
  \cite{Sjogreen} and \cite{BaezaMuletZorio2015}.
    Given $P\in\mathcal{GC}$ and $R\geq r$
we build a set of points $\mathcal{N}(P)=\{N_1,\dots,N_{R+1}\}$ by first selecting the horizontal or vertical direction and considering the first  $R+1$
 interior intersections
of the line determined by $P$ and $N(P)$ with the mesh lines corresponding to 
the selected direction. 
The choice of one direction or the other depends on the slope of the 
normal line, so that the use of interior information is maximized. Figure
 \ref{fig:CxyW} illustrates this approach.
 
Let us remark that in \cite{BaezaMuletZorio2015} the spacing of the 
   nodes at the normal lines used for extrapolation was forced to be
   not smaller than the distance between the ghost cell 
   and the boundary cell that defined the normal line. The reason was that for 
   linear (first order) extrapolation stability issues appear if that constraint 
   is not fulfilled. Numerical experimentation shows that for higher
   order extrapolation there are no stability issues if the nodes 
   (including the boundary condition, if any), are separated between them a distance 
   not smaller than the mesh size, regardless of the position of the ghost cell.
   This results on a more accurate approximation as the node spacing is smaller in this case.
 
We denote by $v=(v_1,v_2)$ the vector determined by $P$ and
$N(P)$, and $\mathcal{S}_q=\{N_{q,1},\dots,N_{q,R+1}\}$ the
closest set of points to $N_q$ from the computational domain sharing
the same coordinate than $N_q$, this depending on the angle of the
normal line as explained in \cite{BaezaMuletZorio2015}. Using the 1D stencil 
$\mathcal{S}_q$ an approximation of the numerical solution at $N_q$ is computed 
 by the same means used for extrapolation at the ghost cell, described 
 in Section \ref{sep}.

Since we are now concerned on performing a weighted
  extrapolation and the computation of
smoothness indicators to build the weights can be very computationally
expensive if the data is not equally spaced, as it will generally happen
on Dirichlet boundaries, we perform an additional step
before extrapolating the value at the ghost cell in order to generate
a new stencil that includes the boundary node and is composed by equally spaced points.

Therefore, if Dirichlet conditions are prescribed, we use the data obtained in
$N_q,$ $1\leq q\leq R+1$ to perform 1D interpolations at the points
$P_q,$ $1\leq q\leq R,$ where $P_q=(P_0^x+qh_x,
P_0^y+q\frac{v_2}{v_1}h_x),$ $0\leq q\leq R$ if $|v_1|\geq|v_2|$ or
$P_q=(P_0^x+q\frac{v_1}{v_2}h_y,P_0^y+qh_y),$ $0\leq q\leq R,$
otherwise, and use the data from the stencil
$\mathcal{S}(P)=\{P_0, P_1, \ldots, P_{R}\}$ to extrapolate it
at the ghost cell $P$.

In case of outflow conditions, we
extrapolate directly the data from the stencil
$\mathcal{S}(P)=\{N_1, N_2, \dots, N_{R+1}\}$. See
    Figure \ref{fig:CxyW} for graphical examples.

If the boundary does not change with time the elements involved in the computation of the extrapolated value at $P$ 
are determined only once at the beginning of the simulation. 

Summarizing, the complete 
extrapolation procedure is done in two stages in the case of outflow 
boundaries and in three in the case of Dirichlet boundaries: in a first
 step data located at the normal lines is computed from
 the numerical solution by (horizontal or vertical) 1D
 interpolation; in the second one, only performed on Dirichlet
 boundaries, the nodes obtained in the first step are used to
 interpolate at new points at the normal line so that they 
 compose an equally spaced stencil together with the point $N(P)$; 
 finally, values for the ghost cell
 are obtained by 1D extrapolation along the normal line
 obtained in the first stage (in case
 of outflow boundary) or the second stage (in case of Dirichlet
 boundary). Note that stencils with equally spaced notes are used in all 
 the above approximations.

We next introduce an extrapolation procedure that takes into account 
the potential presence of discontinuities in the numerical solution.

\section{Extrapolation}
\label{sep}

As stated in \cite{BaezaMuletZorio2015}, special care must be taken
when performing extrapolation at the boundary
mainly due to the potential lack of regularity
in the solution and to the stability restrictions induced by small-cut cells. 
It was also seen there that the classical ENO and WENO methods are not suitable 
for this task
 and a new technique that overcame the above
issue was developed. The method was based on a Boolean criterion for the selection of the nodes that ultimately compose the stencil used to extrapolate. 
A decision on whether a candidate node was added to the stencil or discarded 
was made through an smoothness analysis and controlled via threshold values.
The method produces high order extrapolations when the solution is smooth 
and takes into account the possible presence of discontinuities.

We now present a new technique, which can be considered as an
evolution of the thresholding method, based on the computation of
dimensionless and scale independent weights that outperforms the aforementioned method
and depends on less parameters.

\subsection{Weighted extrapolation}
\label{sep:w}

Consider a stencil of not necessarily
equally nodes $x_0 < \dots < x_r$ and their corresponding nodal values
$u_j=u(x_{j})$. Denote $J=\{0, \dots, r\}$ and $X=\{x_j\}_{j\in J}$ and let $x_*$ be the node 
where we wish to interpolate
and $j_0$ the interior node which is closest to $x_*$, i. e.,
$$j_0=\argmin_{j\in J}|x_j-x_*|.$$
 The goal is to  approximate the value that $x_*$ should
have, based on the information of the ``smoothest'' substencil and the
node $x_{j_0}$. We define
inductively the following set of indexes:
$$J_0=\{j_0\},\textnormal{ and }X_0=\{x_j\}_{j\in J_0}=\{x_{j_0}\}$$

Assume we have defined $J_{k-1}$, then $J_k$ is defined by
$$J_k=J_{k-1}\cup\{j_k\}\textnormal{ and }X_k=\{x_j\}_{j\in J_k}$$
where
$$x_{j_k}=\argmin_{j\in J\setminus J_{k-1}}|x_j- x_*|.$$
That is, we add the closest node to $x_*$ from the remaining nodes to
choose as we increase $k$. $X_k$ and $J_k$ can be defined for $0\leq
k\leq r$.

By construction, it is clear that the set $X_k$
can be written as a sequence of nodes with 
successive indexes, i.e., stencils:
$$X_k=\{x_{i_k+j}\}_{j=0}^k$$
for some $0\leq i_k\leq r-k$, $0\leq k\leq r$.

Now, for each $k$, $0\leq k\leq r$, we define $p_k$ as the interpolating
polynomial of degree at most $k$ such that
$p_k(x_{i_k+j})=u_{i_k+j},$ $\forall j,$ $0\leq j\leq k$.
Given a set of weights $\{\omega_k\}_{k=1}^r$ such that $0\leq\omega_k\leq
1$, we define the following recurrence:

\begin{equation}\label{eq:weightsrec}
\begin{array}{rcl}
  u^{(0)}_*&=&p_0(x_*)=u_{i_0},\\
  u^{(k)}_*&=&(1-\omega_k)u^{(k-1)}_*+\omega_kp_k(x_*),\quad1\leq k\leq r.
  \end{array}
  \end{equation}
  
We define the final result of the weighted extrapolation as
$$u_*:=u^{(r)}_*,$$
which will be taken as an approximation for the value $u(x_*)$.

The idea is to increase the degree of the interpolating polynomial only if the
solution in the corresponding stencil is smooth, 
and therefore the chosen weights should verify that $\omega_k\approx0$ if the
stencil $J_k$ crosses a discontinuity and $\omega_k\approx 1$
if the data from the stencil is smooth. We will show below a weight
construction that verifies that property as well as the capability of
preserving the accuracy order of the extrapolation in case of
smoothness.

From now on, we will assume that the nodes $X$ are equally spaced and
define $h=x_{j+1}-x_j$.

For each $1\leq k\leq r$, we define a slight modification of the
Jiang-Shu smoothness indicator \cite{JiangShu96} associated to the stencil $J_k$ as the
following value:
$$I_k=\frac{1}{r}\sum_{\ell=1}^k\int_{x_0}^{x_r}h^{2\ell-1}p_k^{(\ell)}(x)^2dx.$$

Now, given $\displaystyle1\leq r_0\leq E\left(\frac{r}{2}\right)$,
where $\displaystyle E(x)=\max\mathbb{Z}\cap(-\infty,x]$, we will seek
for a smoothness zone along the stencils of $r_0+1$ points as a
reference.

This procedure will work correctly if there is only one discontinuity in the
stencil, and the restriction $\displaystyle r_0\leq
E\left(\frac{r}{2}\right)$ is set in order to avoid a stencil overlapping,
since a discontinuity might eventually be in the overlapping zone and
thus none of the stencils would include smooth data. Define
$$IS_k=\min_{0\leq j\leq
  r-k}\frac{1}{r}\sum_{\ell=1}^{r_0}\int_{x_0}^{x_r}
h^{2\ell-1}q_{k,j}^{(\ell)}(x)^2dx,\quad1\leq k\leq r_0,$$
where $q_{k,j}$ is the polynomial of degree at most $k$ such that
$q_{k,j}(x_{j+i})=u_{j+i}$ for $0\leq i\leq k,$ $0\leq j\leq r-k$.

There are many possibilities for defining the weights in \eqref{eq:weightsrec}.
We next introduce several possibilities that might be suitable for different 
scenarios. We start with two designs that we name Simple Weights (SW) 
and Improved Weights (IW) that show good performance for problems that have smooth solutions
but might misbehave otherwise. Albeit the interest of such approaches is reduced to academic 
problems we describe them in sections \ref{ap:SW} and \ref{ap:IW} with some detail 
as they are illustrative of the ideas behind other methods described in the rest of the section, 
which are much more oriented to more challenging problems that may include the 
typical non-smooth features of hyperbolic problems.

Numerical experimentation shows that the use of the SW and IW methods 
in complex problems can lead to poor results, probably because of 
the low numerical viscosity introduced by the methods at 
the boundary. For this reason, we
  introduce new weight designs which are derived from the
  two aforementioned extrapolation techniques. These are the unique
  weight (UW),
where the extrapolation is performed by computing one weight, a tuned
version ($\lambda$-UW), where the tuning parameter $\lambda$
can magnetize the weight to 0 or 1 in order to improve the quality of
the extrapolation depending on the context we are working in, and the
global average weight (GAW), also based on a single weight, but more
robust than the UW version.

Finally, since some stability issues may appear if the
extrapolation is solely based on Lagrange extrapolation, we introduce
a least-squares extrapolation procedure, which can be ultimately
combined with any of the the weighted extrapolation techniques. We refer to this
combination as Weighted Least Squares (WLS) so that WLS-X denotes the WLS technique
combined with weights computed by method X. 

\subsection{Simple weights (SW)}
\label{ap:SW}

The weights are defined as follows
\begin{equation}\label{pesos}
  \begin{split}
  &\omega_k=1-\left(1-\left(\frac{IS_k}{I_k}\right)^{s_1}\right)^{s_2},\quad
1\leq k\leq r_0,\\
&\omega_k=\min\left\{1-\left(1-\left(\frac{IS_{r_0}}{I_k}\right)^{s_1}\right)^{s_2},1\right\},\quad
r_0+1\leq k\leq r.
  \end{split}
\end{equation}

A small positive number $\varepsilon>0$ is added to each smoothness indicator in order 
  to avoid the denominator to become zero (in all our experiments, we take
  $\varepsilon=10^{-100}$).
The parameter $s_1$ enforces the convergence to 0 when
the stencil is not smooth, while the parameter $s_2$ enforces the
convergence to 1 when it is smooth.

It can be shown that for a smooth stencil, if there exists some $1\leq
k_0\leq r_0$ such that $|u^{(k_0)}|>>0$ around the stencil, then
$$\omega_k=1-\mathcal{O}(h^{s_2})$$
and if the stencil crosses a discontinuity, then
$$\omega_k=\mathcal{O}(h^{2s_1}).$$

\begin{figure}[htb]
  \centering
  \includegraphics[width=0.5\textwidth]{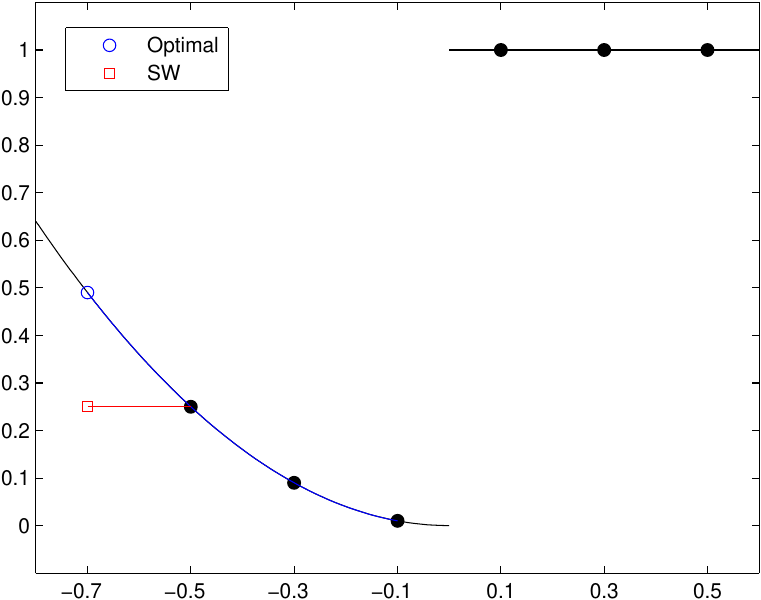}
  \caption{Comparison of the SW extrapolation results against the
    expected optimal result.}
  \label{fig:swissue}
\end{figure}

A drawback of this weight design, apart from the loss of accuracy
when all the
$k$-th derivatives $1\leq k\leq r_0$ vanish near the stencil, is
the fact that sometimes the optimal order
cannot be attained regardless of the values of $s_1$ and
$s_2$. Consider for instance $r=5$ and $r_0=2$ and the function
$$u(x)=\left\{\begin{array}{cc}
    x^2,&x\leq0\\
    1,&x>0
    \end{array}\right.$$
Now we take the nodes $x_i=-0.5+0.2i,$ $0\leq i\leq 5$ and the
corresponding values $u_i$ are $U=\{0.25, 0.09, 0.01, 1, 1,
1\}$. Since one of the substencils of three points contains $\{1, 1,
1\}$, it is clear then that $IS=0$ and therefore $\omega_k=0,$
$\forall k,$ $1\leq k\leq 5$. If we performed a weighted
extrapolation to $x_*=-0.7$ then it would be obtained the
nodal value from the corresponding closest node,
$x_0=-0.5$, that is,
$u_*=u_0=0.25$, while the most reasonable thing to do would be to
perform a
second order extrapolation taking the first three nodes (take all the
nodes from the correct side of the discontinuity), which gives
$0.49$. The aforementioned scenario
  is depicted in Figure \ref{fig:swissue}.

This has occurred because in this case the smoothest substencil belongs
to the other side of the discontinuity, and thus the information about
the derivatives is wrong. We next present an alternative weight design
that overcomes the above issue.

\subsection{Improved weights (IW)}
\label{ap:IW}

The previously mentioned issue can be solved through
  the following modification of the weights design. In this new weight
design we will define and combine additional parameters in order to
ensure that, for stencils with discontinuities, the information is taken
from the correct side of the discontinuity. Let us now consider the
maximum value of the up to $r_0$-th order smoothness indicators:

$$IM_k=\max_{0\leq j\leq
  r-k}\frac{1}{r}\sum_{\ell=1}^{r_0}\int_{x_0}^{x_r}
h^{2\ell-1}q_{k,j}^{(\ell)}(x)^2dx,\quad1\leq k\leq r_0.$$

We define

$$\sigma_k:=\min\left\{\frac{IS_{\min\{k,r_0\}}}{I_k},1\right\},$$
$$\tau_k:=\frac{I_k}{IM_k},$$
and
$$\rho_k=\tau_k\left(\frac{1-\sigma_k}{\sigma_k}\right)^d,\quad
d\geq\frac{r}{2},$$
to finally define our new weights as
\begin{equation}\label{pesosmillorats}
  \omega_k=\frac{1}{1+\rho_k}.
\end{equation}

We now see in detail the reason for this choice. These are the possible
values in terms of powers of $\mathcal{O}(h)$ that both quotients
between smoothness indicators can take:

$$\sigma_k=\left\{\begin{array}{rl}
    1-\mathcal{O}(h) & \textnormal{if there is smoothness,} \\
    \mathcal{O}(h^2) & \textnormal{if } X_k \textnormal{ crosses a
      discontinuity.}
    \end{array}\right.$$
However, if there is not smoothness but still $X_k$ does not cross a
discontinuity there may be two possibilities: $IS_{\min\{k,r_0\}}$ is
obtained in the ``correct'' side or in the ``wrong'' one. In the first
case, we would have
$\displaystyle\frac{IS_{\min\{k,r_0\}}}{I_k}=1+\mathcal{O}(h)$, as desired, but
otherwise, if the derivatives are different from one side and another
in a random fashion, that would lead to quotients with values in a
random fashion as well. Here is thus the importance of the second
quotient in order to fix such issue:
$$\tau_k=\left\{\begin{array}{rl}
    1+\mathcal{O}(h) & \textnormal{if there is smoothness,} \\
    \mathcal{O}(h^2) & \textnormal{if there is not smoothness and }
    X_k \textnormal{ belongs to a smooth zone,} \\
    \mathcal{O}(1) & \textnormal{if there is not smoothness and } X_k
    \textnormal{ crosses a discontinuity,}
    \end{array}\right.$$
where $\mathcal{O}(1)$ means that it is a value comprised between 0
and 1 in a random fashion, but not $\mathcal{O}(h)$ since in that case
it is a quotient between two smoothness indicators, both in non-smooth
zones. Note from that the term $IM_k$ from the definition of $\tau_k$
can be replaced by
$I_r$, since in this case $\tau$ will still verify the above
properties and less smoothness indicators will be required to be
computed. In practice, we will work through this modification.

One can show that this way

$$\rho_k=\left\{\begin{array}{rl}
    \mathcal{O}(h^d) & \textnormal{if } X_k \textnormal{ belongs to a
      smooth zone,} \\
    \mathcal{O}(h^{-2d}) & \textnormal{if } X_k \textnormal{ crosses a
      discontinuity,} \\
    \end{array}\right.$$

as desired, and thus

$$\omega_k=\left\{\begin{array}{rl}
    1-\mathcal{O}(h^{d}) & \textnormal{if } X_k \textnormal{ belongs to a
      smooth zone,} \\
    \mathcal{O}(h^{2d}) & \textnormal{if } X_k \textnormal{ crosses a
      discontinuity,} \\
    \end{array}\right.$$
except when all derivatives up to the $r_0$-th one are close to be
zero in some of the two discontinuity sides, where the value of the
weight cannot be predicted, although it will take values more likely
close to 0 rather than 1.

We can fix this issue by redefining $\sigma_k$ as
\begin{equation}\label{beta}
  \sigma_k=\frac{IS_k+\beta}{I_k+\beta},
\end{equation}
where $\beta$ is assumed to be a small quantity,
$\beta=\mathcal{O}(h^b),$ $b\leq 2r_0$. We propose two different
choices of $\beta$ satisfying that condition:
\begin{itemize}
  \item $\beta=\lambda^2 h^{2r_0}$, where $\lambda$ is a parameter
    proportional to the scaling of the solution. That is, if we
    re-scale the data by a factor of $\mu$, then the value
    $\lambda$ should be replaced by $\mu\lambda$.
  \item In our context, we have knowledge about a wide range
    of point values of a function (the numerical data from a
    computational domain in a certain scheme). Hence, we can naturally
    define a quantity $\mathcal{O}(h^{2r_0})$ depending directly on the
    data of the problem.

    For instance, if we assume we have a 1D simulation with data $U_j,$
    $0\leq j\leq n$, then one can define
    $$\beta=\left(\frac{1}{n}\sum_{j=1}^n|u_j-u_{j-1}|\right)^{2r_0}
    =TV(u)^{2r_0}h^{2r_0}.$$
    The above value can be generalized to any dimension as a global
    average of all the
    absolute values of all the directional undivided differences (in
    all directions).
    
    When discontinuities and zeros at the first derivative are
    supposed to be in a region of measure 0, the above value verifies
    $\beta=\mathcal{O}(h^{2r_0})$, while keeping
    $\sigma_k$ scaling independent. This argument is also valid even
    when the first derivative is zero almost everywhere, but
    there is a discontinuity on the data, which is also a common case
    as initial condition for shock problems. Another valid case is
    when, despite having zeros in the derivative in a non-null region,
    there is a non-null region having non-zero derivatives as well
    provided that the discontinuities, if any, are located in a null
    region.

    If one wants to keep some stricter control of the above parameter
    due to the presence of very strong discontinuities which might
    make $\beta$ too big, one can always consider a tuning parameter
    $\kappa$ (in this case independent of the scaling as well) and
    redefine $\sigma_k$ as
    $$\sigma_k=\frac{IS_k+\kappa\beta}{I_k+\kappa\beta}.$$
\end{itemize}
If one uses the above technique to avoid a loss of accuracy near zeros
on the corresponding derivatives and makes sure that the parameters,
if any, are well tuned, it makes no sense to use a smoothness control
stencil longer than a two-points one, and thus for that case we will
always use $r_0=1$. Moreover, in this particular case it is no longer
needed to define $\sigma_k$ such that does not surpass the unity,
since we have the following result.

\begin{propr}
  If $r_0=1$ then
  $$\alpha_k:=\frac{I_1}{I_k}$$
  verifies $0\leq\alpha_k\leq1,$ $1\leq k\leq r$.
\end{propr}
\begin{proof}
Given $f,g\in L^2([x_{i-1},x_i])$ we define the
following scalar products
$$\langle f,g\rangle_i=\int_{x_{i-1}}^{x_i}f(x)g(x)dx$$
and their induced norms as
$$\|f\|_i^2=\langle f,f\rangle_i,\quad1\leq i\leq r.$$
Now, taking into account that $x_i-x_{i-1}=h$, we have
\begin{equation*}
\begin{split}
I_k&=\frac{1}{r}\sum_{\ell=1}^k\int_{x_0}^{x_r}h^{2\ell-1}p_k^{(\ell)}(x)^2dx\geq
\frac{1}{r}\int_{x_0}^{x_r}hp_k'(x)^2dx=
\frac{1}{r}\sum_{i=1}^rh\int_{x_{i-1}}^{x_i}p_k'(x)^2dx\\
&=\frac{1}{r}\sum_{i=1}^r
\int_{x_{i-1}}^{x_i}dx\int_{x_{i-1}}^{x_i}p_k'(x)^2dx=
\frac{1}{r}\sum_{i=1}^r\|f_i\|_i^2\|g_i\|_i^2,
\end{split}
\end{equation*}
where $f_i(x)=1$ and $g_i(x)=p_k'(x)$ for $x\in[x_{i-1},x_i]$. By the
Cauchy-Schwarz Inequality
$$\langle f_i,g_i\rangle_i^2\leq\|f_i\|_i^2\|g_i\|_i^2$$
we have
\begin{equation*}
\begin{split}
I_k&\geq\frac{1}{r}\sum_{i=1}^r\|f_i\|_i^2\|g_i\|_i^2\geq
\frac{1}{r}\sum_{i=1}^r\langle f_i,g_i\rangle_i^2=
\frac{1}{r}\sum_{i=1}^r\left(\int_{x_{i-1}}^{x_i}f_i(x)g_i(x)dx\right)^2\\
&=\frac{1}{r}\sum_{i=1}^r\left(\int_{x_{i-1}}^{x_i}p'_k(x)dx\right)^2
=\frac{1}{r}\sum_{i=1}^r(p_k(x_i)-p_k(x_{i-1}))^2=
\frac{1}{r}\sum_{i=1}^r(u_i-u_{i-1})^2\\
&\geq\frac{1}{r}\sum_{i=1}^r\min_{1\leq j\leq r}(u_j-u_{j-1})^2=
\min_{1\leq j\leq r}(u_j-u_{j-1})^2=I_1.
\end{split}
\end{equation*}
\qed
\end{proof}

With these modifications, weighted extrapolation will not suffer from
a loss of accuracy in any case and it will have the optimal order in
presence of sharp discontinuities.

If we now apply this technique, it will still capture well sharp
discontinuities while keeping the highest order as possible when there
is a discontinuity in the stencil.

\subsection{Examples}
\begin{figure}[htb]
  \centering
  \includegraphics[width=0.5\textwidth]{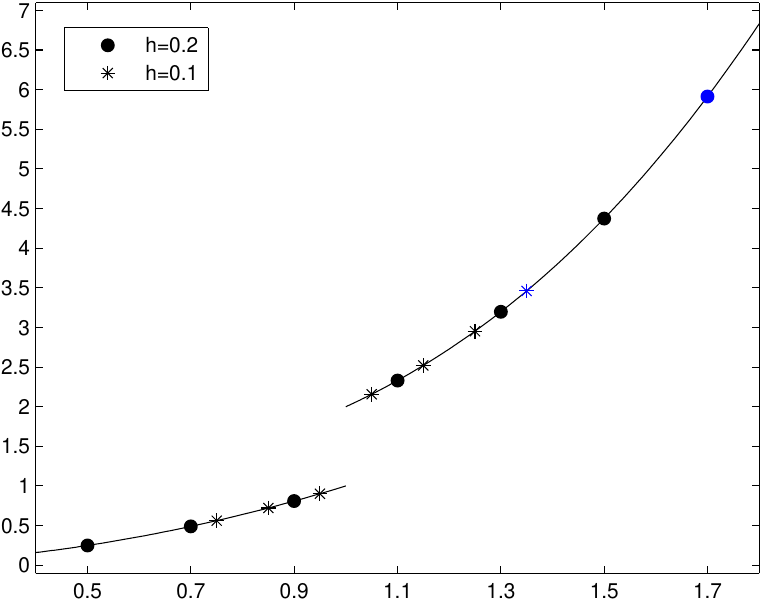}
  \caption{Illustration of the grid configuration in Example 1.}
  \label{fig:gex1}
\end{figure}

\textbf{Example 1.} Let $u:\mathbb{R}\rightarrow\mathbb{R}$ a function
defined by
$$u(x)=\left\{\begin{array}{ccc}
    x^2 & \textnormal{if} & x\leq 1, \\
    1+x^3 & \textnormal{if} & x>1.
\end{array}\right.$$
We study numerically the accuracy behavior of our scheme in a six
points stencil around the discontinuity point $x=1$ when
$h\to0$. Given $h>0$ we select as stencil the set of nodes
$x_i=1+(-2.5+i)h,$ $0\leq i\leq 5$. Figure
  \ref{fig:gex1} shows a graphical example of the grid points for
  $h=0.2$ and $h=0.1$.

We start with $h_0=0.04$ and define $h_i=h_{i-1}/2$ for $1\leq i\leq
6$ and compute the exact errors extrapolating at $x_*=1+3.5h$ using the
above techniques. A successful one should give third order accuracy
since there are available three points from the right side of the
discontinuity.

The results obtained with the simple weights (\ref{pesos}) for
$s_1=s_2=3$ and $r_0=2$ are shown in Table \ref{taulapesos}.

\begin{table}[htb]
  \centering
  \begin{tabular}{|c|c|c|}
    \hline
    Resolution & Error & Order \\
    \hline
    $h_0$ & 5.72E$-$2 & $-$ \\
    \hline
    $h_1$ & 1.18E$-$2 & $1.48$ \\
    \hline
    $h_2$ & 7.40E$-$3 & $1.32$ \\
    \hline
    $h_3$ & 3.24E$-$3 & $1.19$ \\
    \hline
    $h_4$ & 1.51E$-$3 & $1.10$ \\
    \hline
    $h_5$ & 7.26E$-$4 & $1.06$ \\
    \hline
    $h_6$ & 3.56E$-$4 & $1.03$ \\
    \hline
  \end{tabular}
  \caption{Simple weights (\ref{pesos}), example 1.}
  \label{taulapesos}
\end{table}
It can be clearly seen that the optimal accuracy is not attained by
this weight design.

The technical reason for this happening is that the derivative of $u$
for $x\neq1$ is
$$u'(x)=\left\{\begin{array}{ccc}
    2x & \textnormal{if} & x<1, \\
    3x^2 & \textnormal{if} & x>1,
\end{array}\right.$$
and thus $\lim_{x\to 1^-}u'(x)=2$ and $\lim_{x\to 1^+}u'(x)=3$, hence
since the ``smoothest'' information is taken from the left side of the
discontinuity but the actual information should be taken from the
right side of the discontinuity, where there is smoothness as well.
Therefore, the weights $\omega_k,$ $1\leq k\leq 2$, converge to
$(1-(\frac{2}{3})^3)^3$ as $h\to0$ rather than 1 for the above
explained reason.

Now, we repeat the same test using the new weights defined in
(\ref{pesosmillorats}) for $d=3$ and $r_0=1$ and we present in Table
\ref{taulapesosmillorats} the errors for the same setup as above,
where it can be clearly seen that the optimal third order accuracy is
obtained.

\begin{table}[htb]
  \centering
  \begin{tabular}{|c|c|c|}
    \hline
    Resolution & Error & Order \\
    \hline
    $h_0$ & 3.89E$-$4 & $-$ \\
    \hline
    $h_1$ & 4.82E$-$5 & $3.01$ \\
    \hline
    $h_2$ & 6.01E$-$6 & $3.00$ \\
    \hline
    $h_3$ & 7.50E$-$7 & $3.00$ \\
    \hline
    $h_4$ & 9.38E$-$8 & $3.00$ \\
    \hline
    $h_5$ & 1.17E$-$8 & $3.00$ \\
    \hline
    $h_6$ & 1.46E$-$9 & $3.00$ \\
    \hline
  \end{tabular}
  \caption{Improved weights (\ref{pesosmillorats}), example 1.}
  \label{taulapesosmillorats}
\end{table}

\begin{figure}[htb]
  \centering
  \includegraphics[width=0.5\textwidth]{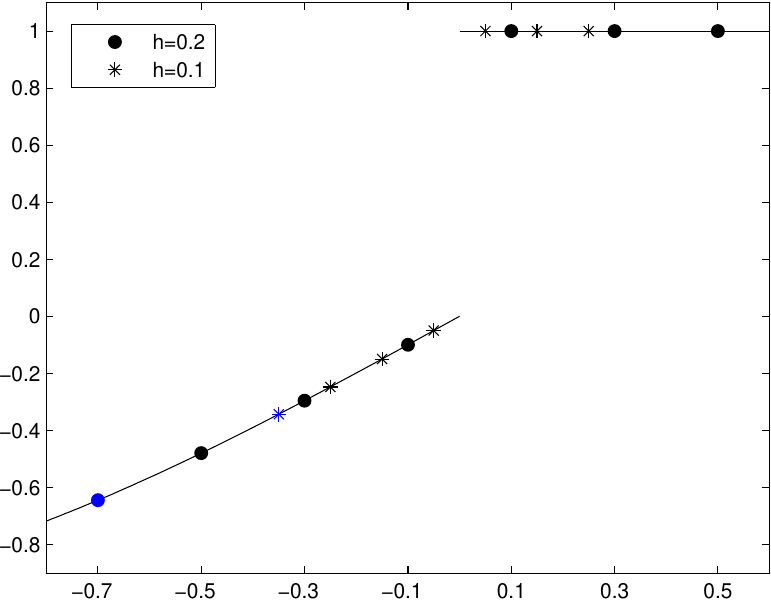}
  \caption{Illustration of the grid configuration in Example 2.}
  \label{fig:gex2}
\end{figure}

\textbf{Example 2.} We now consider an example where one derivative
vanishes. This example is very similar to the one presented to motivate the
definition of the improved weights. In this case, we define
$u:\mathbb{R}\rightarrow\mathbb{R}$ as
$$u(x)=\left\{\begin{array}{ccc}
    \sin(x) & \textnormal{if} & x\leq 0, \\
    1 & \textnormal{if} & x>0.
\end{array}\right.$$
We define now the grid points as $x_i=(-2.5+i)h,$
  $0\leq i\leq 5$, whose configuration is depicted in Figure
   \ref{fig:gex2} for $h=0.2$ and $h=0.1$.

Taking the same values
of $h$ as in the above experiment, we now extrapolate at $x_*=-3.5h$
and we obtain the results in Table \ref{taulapesosmilloratsnobeta}
 using the improved weights with the same parameters as above.
This order decay is now due to the fact that one of the derivatives
vanishes.

\begin{table}[htb]
  \centering
  \begin{tabular}{|c|c|c|}
    \hline
    Resolution & Error & Order \\
    \hline
    $h_0$ & 3.97E$-$2 & $-$ \\
    \hline
    $h_1$ & 2.00E$-$2 & $0.99$ \\
    \hline
    $h_2$ & 1.00E$-$2 & $1.00$ \\
    \hline
    $h_3$ & 5.00E$-$3 & $1.00$ \\
    \hline
    $h_4$ & 2.50E$-$3 & $1.00$ \\
    \hline
    $h_5$ & 1.25E$-$3 & $1.00$ \\
    \hline
    $h_6$ & 6.25E$-$4 & $1.00$ \\
    \hline
  \end{tabular}
  \caption{Improved weights (\ref{pesosmillorats}), example 2.}
  \label{taulapesosmilloratsnobeta}
\end{table}

To fix this, we use the weight modification suggested in (\ref{beta})
by taking $\beta=h^2$, obtaining third order accuracy as can be observed in 
table \ref{taulapesosmilloratsbeta}.

\begin{table}[htb]
  \centering
  \begin{tabular}{|c|c|c|}
    \hline
    Resolution & Error & Order \\
    \hline
    $h_0$ & 6.38E$-$4 & $-$ \\
    \hline
    $h_1$ & 7.99E$-$5 & $3.00$ \\
    \hline
    $h_2$ & 1.00E$-$5 & $3.00$ \\
    \hline
    $h_3$ & 1.25E$-$6 & $3.00$ \\
    \hline
    $h_4$ & 1.56E$-$7 & $3.00$ \\
    \hline
    $h_5$ & 1.95E$-$8 & $3.00$ \\
    \hline
    $h_6$ & 2.44E$-$9 & $3.00$ \\
    \hline
  \end{tabular}
  \caption{Improved weights (\ref{pesosmillorats}), $\beta=h^2$, example 2.}
  \label{taulapesosmilloratsbeta}
\end{table}

\subsection{Unique weight extrapolation (UW)}

In this section we propose a method that uses only one weight to decide 
the extrapolation method attending to the global smoothness in the stencil $X$ that contains all $r+1$ nodes. 
The idea is that switching to a low order reconstruction as soon as a lack of smoothness is detected
increases the robustness of the procedure in the non-smooth case
when the extrapolation is performed in the context of
  PDEs, where discontinuities get smeared,
while maintaining high order in the smooth case, besides being
  more computationally efficient.

  As our purpose, apart from achieving robustness, is
    also designing an efficient method, we will propose 
    for the UW method a simplified version of the smoothness indicators 
    introduced previously. We will also
    introduce in Section \ref{sec:tunedw} an additional (and optional) tuning parameter, with which we can map
the original weight $\omega$ to another one through a transformation
that magnetizes values
far from 1 (but still far from 0 as well) to 0. This variant is particularly useful in problems
with strong shocks --as in this case the weights should be very
close to 0-- or in problems with complex smooth structure --where the results are better if the weights
are close to 1 near them--.

The computation of the weights in the previous cases
  implies the use
of logical structures since a minimum has to be computed. We overcome
this issue as well by designing a weight capable of capturing well
the discontinuities while keeping high order accuracy on smooth
zones. This will be discussed in section \ref{sec:gaw}.

The procedure for the extrapolation with only one weight is performed
in the following sense: Instead of
gradually increasing the degree of the interpolating polynomials, we
will just average the constant extrapolation ($k=0$) and maximum
degree extrapolation ($k=r$), that is, we will consider
$$u_*=(1-\omega)p_0(x_*)+\omega p_r(x_*)=(1-\omega)u_{i_0}+\omega
p_r(x_*),$$
where
$$\omega=\min\left\{1-\left(1-\left(\frac{IS_{r_0}}{I_r}\right)^{s_1}\right)^{s_2},1\right\}.$$

It can be shown that for a smooth stencil, if there exists some $1\leq
k_0\leq r_0$ such that $|u^{(k_0)}|>>0$ around the stencil, then
$$\omega=1-\mathcal{O}(h^{s_2})$$
and if the stencil crosses a discontinuity, then
$$\omega=\mathcal{O}(h^{2s_1}).$$

In order to lower the computational cost and ensure that
$0\leq \omega\leq 1$ without having to artificially bound it by 1 when
$r_0>1$, we can replace the
definition of $I_r$, which is a smoothness indicator of the whole $r+1$
points stencil, by the average of all smoothness indicators of the
substencils of $r_0+1$ points, i.e.:

\begin{equation}\label{newsi}
  I_r^*:=\frac{1}{r-r_0+1}\sum_{j=0}^{r-r_0}I_{r_0,j},$$
where
$$I_{r_0,j}=\frac{1}{r_0}\sum_{\ell=1}^{r_0}\int_{x_j}^{x_{r_0+j}}h^{2\ell-1}q_{r_0,j}^{(\ell)}(x)^2dx.
\end{equation}

Then one can define

$$\omega=1-\left(1-\left(\frac{IS_{r_0}}{I_r^*}\right)^{s_1}\right)^{s_2},$$
which in this case it clearly verifies $0\leq\omega\leq1$.

Under the hypothesis $\exists k_0\in\mathbb{N}$, $1\leq
k_0\leq r_0$ such that $|u^{(k_0)}|>>0$ around the stencil, then
$$u_*=u(x^*)+\mathcal{O}(h^{r'+1}),$$
where $r'=\min\{s_2(r_0-k_0+1), r\}$.

\subsubsection{Tuned weight ($\lambda$-UW)}\label{sec:tunedw}

When the stencil
includes globally some smeared discontinuity, one can map a $\omega$
value such that $\omega>>0$ and $\omega<<1$ to a new one
$\widetilde{\omega}\approx0$, as desired in this case, but that at same
time $\widetilde{\omega}\approx1$ when $\omega\approx1$ as well.

Let $\tw=F_{\lambda}(\omega)$ be
$$\tw:=\left\{\begin{array}{ccc}
  \frac{e^{\lambda\omega}-1}{e^{\lambda}-1} & \textnormal{if} &
  \lambda\neq0\\
  \omega & \textnormal{if} & \lambda=0
  \end{array}\right.,$$
for some $\lambda\in\mathbb{R}$. It can be proven that
$\forall\omega\in[0,1],$ 
$G_\omega\in\mathcal{C}^{\infty}(\mathbb{R})$, where
$G_\omega(\lambda):=F_{\lambda}(\omega)$. The larger $\lambda$ is, the
stricter the discontinuity detection filter will be. The lower
(negative) $\lambda$ is, the more permissive it will
be.

Assume smoothness conditions, therefore $\omega=1-\delta$, where
$\delta=\mathcal{O}(h^{s_1})$. Assume $s_1\geq r$ in order to
achieve the maximum order accuracy. We have:
\begin{equation*}
\begin{split}
e^{\lambda\omega}-1&=e^{\lambda(1-\delta)}-1=e^{\lambda}e^{-\lambda\delta}-1=
e^{\lambda}e^{-\lambda\delta}-1
=e^{\lambda}e^{-\lambda\delta}-e^{\lambda}+e^{\lambda}-1 \\
&=(e^{\lambda}-1)+e^{\lambda}(e^{-\lambda\delta}-1).
\end{split}
\end{equation*}
Therefore
$$\tw=\frac{(e^{\lambda}-1)+e^{\lambda}(e^{-\lambda\delta}-1)}{e^{\lambda}-1}
=1+\frac{e^{\lambda}}{e^{\lambda}-1}(e^{-\lambda\delta}-1).$$
On the other hand, using the Taylor expansion:
$$e^{-\lambda\delta}-1=\sum_{k=1}^{\infty}\frac{(-\lambda\delta)^k}{k!}
=-\lambda\mathcal{O}(h^{s_1}).$$
Therefore
$$\tw=1-\frac{e^{\lambda}}{e^{\lambda}-1}\lambda\mathcal{O}(\Delta
x^{s_1})$$
and the desired accuracy is attained provided that
$0<<\lambda=\mathcal{O}(1)$.

\begin{figure}
\centering
\includegraphics[scale=0.5]{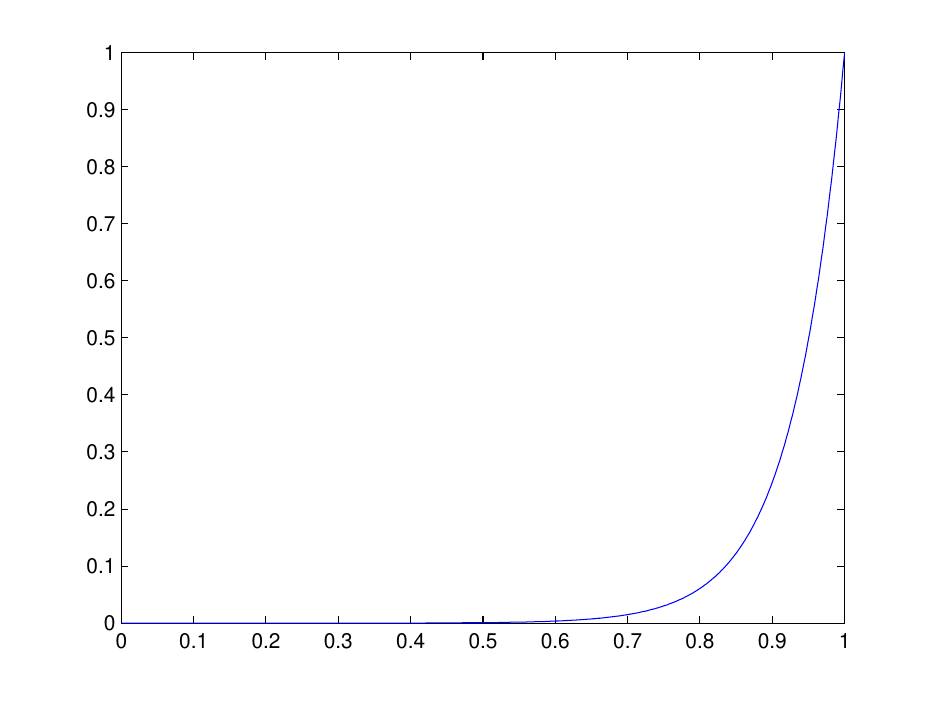}
\caption{Plot of the mapping $\omega\rightarrow\tw$ for $\lambda=14$.}
\label{w}
\end{figure}

Figure \ref{w} shows a plot of the mapping $\omega\rightarrow\tw$ for
$\lambda=14$.

\subsubsection{Global average weight (GAW)}\label{sec:gaw}


Using the weight defined through the smoothness indicator replacement
in (\ref{newsi}), a natural substitute of $IS_{r_0}$ using every available $r_0$-th order smoothness 
indicator is their harmonic mean:
$$IS_{r_0}^*:=\frac{1}{\frac{1}{r-r_0+1}\sum_{j=0}^{r-r_0}\frac{1}{I_{r_0,j}}}.$$
We can thus define the new weight $\omega$ as
$$\omega:=(1-(1-\rho)^{s_1})^{s_2},$$
where
$$\rho:=\frac{IS_{r_0}^*}{I_r^*}=
\frac{\frac{1}{\frac{1}{r-r_0+1}\sum_{j=0}^{r-r_0}\frac{1}{I_{r_0,j}^m}}}
{\frac{1}{r-r_0+1}\sum_{j=0}^{r-r_0}I_{r_0,j}^m}=
\frac{(r-r_0+1)^2}{\left(\sum_{j=0}^{r-r_0}I_{r_0,j}^m\right)
  \left(\sum_{j=0}^{r-r_0}\frac{1}{I_{r_0,j}^m}\right)},$$
and $m$ is a parameter that enforces the convergence to 0 of the
weight in a discontinuity. We next show that $0\leq\rho\leq1$ (and therefore $\omega$
verifies this property as well) as well as the desired properties
both in smooth and non-smooth cases.
\begin{propr}
  $0\leq\rho\leq 1$ and verifies
  $$\rho=\left\{\begin{array}{rl}
    1-\mathcal{O}(\Delta x^{2c_s}) & \textnormal{if the stencil is
    }C^{r_0}\textnormal{ with a }s\textnormal{-th order zero derivative}, \\
    \mathcal{O}(\Delta x^{2ms}) & \textnormal{if the stencil contains a
    discontinuity,}
  \end{array}\right.$$
  with $c_s:=\max\{r_0-s,0\}$. Therefore
  $$\omega=\left\{\begin{array}{rl}
    1-\mathcal{O}(\Delta x^{2s_1c_s}) & \textnormal{if the stencil is
    }C^{r_0}\textnormal{ with a }s\textnormal{-th order zero derivative}, \\
    \mathcal{O}(\Delta x^{2s_2ms}) & \textnormal{if the stencil contains a
    discontinuity.}
  \end{array}\right.$$
\end{propr}
\begin{proof}
  Let $a_j>0,$ $1\leq j\leq k$. We show that the quotient between
  their harmonic mean and their mean is bounded by 1, that is
  $$0\leq\rho=\frac{\frac{1}{\frac{1}{k}\sum_{j=1}^k\frac{1}{a_j}}}{\frac{1}{k}\sum_{j=1}^ka_j}\leq
  1.$$
  Since
  $$\rho=\frac{\frac{1}{\frac{1}{k}\sum_{j=1}^k\frac{1}{a_j}}}{\frac{1}{k}\sum_{j=1}^ka_j}=
  \frac{k^2}{(\sum_{j=1}^k\frac{1}{a_j})(\sum_{j=1}^ka_j)}$$
  it suffices to show that
  $$A:=\left(\sum_{j=1}^k\frac{1}{a_j}\right)\left(\sum_{j=1}^ka_j\right)\geq
  k^2.$$
  Indeed,

  \begin{equation*}
    \begin{split}
      A&=\sum_{i=1}^k\sum_{j=1}^k\frac{a_i}{a_j}
      =k+\sum_{i=1}^k\sum_{j=1}^{i-1}\left(\frac{a_i}{a_j}+\frac{a_j}{a_i}\right)
      =k+\sum_{i=1}^k\sum_{j=1}^{i-1}\frac{(a_i-a_j)^2+2a_ia_j}{a_ia_j}\\
      &=k+2\sum_{i=1}^k(i-1)+\sum_{i=1}^k\sum_{j=1}^{i-1}\frac{(a_i-a_j)^2}{a_ia_j}
      =k^2+\sum_{i=1}^k\sum_{j=1}^{i-1}\frac{(a_i-a_j)^2}{a_ia_j}\geq k^2.
    \end{split}
  \end{equation*}
  Let us now assume that the stencil is $C^{r_0}$ smooth with a $s$-th
  zero in the derivative, then $$I_{r_0,j}=\Delta
  x^{2s}(1+\mathcal{O}(\Delta x^{\max\{r-s,0\}}))=\Delta
  x^{2s}C(1+\mathcal{O}(\Delta x^{c_s})).$$
  Since we have actually proven that
  $$A=k^2+\sum_{i=1}^k\sum_{j=1}^{i-1}\frac{(a_i-a_j)^2}{a_ia_j},$$
  replacing the $a$ terms with the corresponding smoothness indicators
  to the power of $m$ in case of smoothness
  $$I_{r_0,j}^m=\Delta x^{2ms}C^m(1+\mathcal{O}(\Delta x^{c_s}))^m
  =\Delta x^{2ms}C^m(1+\mathcal{O}(\Delta x^{c_s})),$$
  we have, denoting $k:=r-r_0+1$,
  \begin{equation*}
  \begin{split}
    A&=k^2+\sum_{i=1}^k\sum_{j=1}^{i-1}\frac{(\Delta
  x^{2ms}C^m(1+\mathcal{O}(\Delta x^{c_s}))-\Delta
  x^{2ms}C^m(1+\mathcal{O}(\Delta x^{c_s})))^2}{\Delta
  x^{4ms}C^{2m}(1+\mathcal{O}(\Delta x^{c_s}))^2}\\
&=k^2+\sum_{i=1}^k\sum_{j=1}^{i-1}\frac{\Delta x^{4ms}C^{2m}(\mathcal{O}(\Delta x^{c_s}))^2}{\Delta
  x^{4ms}C^{2m}(1+\mathcal{O}(\Delta x^{c_s}))}
=k^2+\sum_{i=1}^k\sum_{j=1}^{i-1}\frac{\mathcal{O}(\Delta
  x^{2c_s})}{1+\mathcal{O}(\Delta x^{c_s})}\\
&=k^2+\mathcal{O}(\Delta x^{2c_s}).
\end{split}
  \end{equation*}
  Therefore,
  $$\rho=\frac{k^2}{A}=\frac{k^2}{k^2+\mathcal{O}(\Delta
    x^{2c_s})}=1+\mathcal{O}(\Delta x^{2c_s}).$$

  Finally, if there is smoothness at least in one substencil with
  corresponding smoothness indicator $I_{r_0, i_0}=\Delta
  x^{2s}C(1+\mathcal{O}(\Delta x^{c_s}))=\mathcal{O}(\Delta x^{2s})$ and a
  discontinuity crosses the stencil then there will be another
  substencil such that its smoothness indicator verifies
  $I_{r_0,j_0}=\mathcal{O}(1)$. Then, the corresponding term in the above
  sum (swapping the indexes if necessary, that is, if $i_0<j_0$):  
  
  $$\frac{(I_{r_0,i_0}^m-I_{r_0,j_0}^m)^2}{I_{r_0,i_0}^mI_{r_0,j_0}^m}=
  \frac{(\mathcal{O}(\Delta x^{2ms})-\mathcal{O}(1))^2}
  {\mathcal{O}(\Delta x^{2ms})\mathcal{O}(1)}
  =\frac{\mathcal{O}(1)}{\mathcal{O}(\Delta
    x^{2ms})}=\mathcal{O}(\Delta x^{-2ms}),$$
  and thus
  $$A=\mathcal{O}(\Delta x^{-2ms}).$$
  Therefore,
  $$\rho=\frac{k^2}{A}=\frac{k^2}{\mathcal{O}(\Delta
    x^{-2ms})}=\mathcal{O}(\Delta x^{2ms}).$$
  \qed
\end{proof}
It follows from the result that to obtain an order of accuracy as large
as possible, and if $r=2r_0$, as it will be our usual choice of $r_0$
from now on, we must take $s_1\geq r_0$ (to prevent from a possible
extreme case $s=r_0-1$). Note that in the worst case scenario, $s\geq r_0$, implies an
unavoidable downgrade of the accuracy as it happens with the original
WENO-JS weights for reconstructions. The only way to overcome this is
issue is to add to each smoothness indicator the $\beta$ defined
above, but it will not be done in the forthcoming tests as we do not
do it either in the definition of the WENO-JS weights in our numerical
solver.

\subsection{Weighted least squares extrapolation (WLS)}

After some 2D order accuracy tests for smooth solutions with different
spacing setups, it has been noticed that depending
on the complexity of the domain, not only high order may not be
achieved using straight Lagrange extrapolation at the ghost cells, but
also the scheme might turn mildly unstable in some extreme cases, this
independently of the choice for the spacing of the normal line points
to be extrapolated at the ghost cell.

To avoid this phenomena, it is necessary to find an alternative to
Lagrange (and weighted) extrapolation that stabilizes the scheme 
while keeping high order accuracy in terms of the global error.
A possibility is to use least squares fitting as described in
\cite{TanShu}, \cite{TanWangShu}. 

Let $R\geq r$ and let $\{(x_i,u_i)\}_{i=0}^R$ be the stencil with nodal data,
where $x_{i+1}-x_i=h$. Let us assume that we want to extrapolate or
interpolate at the
point $x_*$ with a certain $r$-th degree polynomial
$$p(x)=\sum_{j=0}^ra_jx^j$$
using the data from the whole stencil. Since in general there is 
no polynomial of $r$ passing through the $R$ points, we find the polynomial 
$p(x)$ of
degree $r$ which minimizes the error with respect to $\{(x_i,u_i)\}_{i=0}^R$ in the
$L^2$-norm, i. e., we solve $p(x_i)=u_i,$ $0\leq i\leq R$ by least squares,
which will be still a $(r+1)$-th order accurate
approximation if the data is smooth. 

We can now combine the least squares extrapolation, conceived for
smooth regions with the already described techniques 
techniques based on weights. Let $z_*$ be the result of the
least squares extrapolation and $v_*$ the result of some chosen
modality of the weighted extrapolation techniques. 
To this aim we define a weight $\omega$ by using the information of the whole
stencil of $R+1$ points through the $\lambda$-UW or GAW
technique and then take the final result of the extrapolation as
\begin{equation}\label{eq:wls1}
u_*:=\omega z_*+(1-\omega)v_*.
\end{equation}

The simplest and most robust case for the election of $v_*$ is 
to take $v_*=u_{i_0}$. This is ultimately the technique that will be used both for smooth and
non-smooth problems in our numerical experiments, so that we will indicate 
by WLS-X the weighted least squares method with $\omega$ in \eqref{eq:wls1} computed through
the method X.

We next show a step-by-step algorithm of the chosen
  WLS-GAW technique, for a stencil $\{x_i\}_{i=0}^R$ with nodal values
  $\{u_i\}_{i=0}^R$, to be extrapolated at the point $x_*$, is thus:

\begin{enumerate}
  \item Find the index $0\leq i_0\leq R$ such that
      $x_{i_0}$ is the
    closest point to $x_*$. Then $u_{i_0}$ is the reference value.
  \item Compute the least-square polynomial of degree $r$ at $x_*$,
    $p(x_*)$, using the whole stencil data.
  \item Obtain the corresponding smoothness indicators from the
    substencils of size $r_0$.
  \item Compute the global average weight, $w$ from the previously
    computed smoothness indicators.
  \item The final extrapolated value, $u_*$, is then obtained as a
    weighted average of the least-square polynomial and the reference
    value, namely:
    $$u_*=\omega p(x_*)+(1-\omega)u_{i_0}.$$
\end{enumerate}

We next summarize all parameters related to the size of the 
various stencils involved in the extrapolation:

\begin{itemize}
  \item $r+1$: Stencil size for Lagrange
      extrapolation. The accuracy
    order is thus $r+1$ in case of smoothness.
  \item $r_0+1$: Substencil size used in the computation of the
    smoothness indicators.
  \item $R+1$: Stencil size for least-squares extrapolation. In this
    context, we refer to $r$ as the degree of the computed polynomial
    (hence the accuracy order is $r+1$ as well).
\end{itemize}

 Table \ref{lpr} shows all the parameters
 involved in the different weight designs and indications on whether a parameter is involved or not 
 on a particular method.

  \begin{table}[htb]
  \centering
  \begin{tabular}{|l|c|c|c|c|}
    \hline
    Method (columns) / Parameter (rows) & SW & IW & UW & GAW \\
    \hline
    $s_1$: $\omega=\mathcal{O}(h^{\xi s_1})$ if discontinuity & Y & N & Y & Y \\
    \hline
    $s_2$: $\omega=1-\mathcal{O}(h^{\xi s_2})$ if smoothness & Y & N & Y & Y
    \\
    \hline
    $d$:
    \begin{tabular}{l}
      $\omega=\mathcal{O}(h^{\xi d})$ if discontinuity \\
      $\omega=1-\mathcal{O}(h^{\xi d})$ if smoothness
    \end{tabular} & N & Y & N & N \\
    \hline
    $\lambda$:
    \begin{tabular}{l}
      $\omega\to0$ if $\lambda\to+\infty$ \\
      $\omega\to1$ if $\lambda\to-\infty$
    \end{tabular} & N & N & O & N \\
    \hline
    $m$: $\omega=\mathcal{O}(h^{\xi m})$ if discontinuity & N & N & N
    & Y \\
    \hline
  \end{tabular}
  \caption{List of parameters for the weights (Y: yes, N: no, O:
    optional). In each case, $\xi$ is a parameter that depends on the
    method, the substencils size $r_0$, the number of consecutive
    zero-derivatives and other parameters.}
    \label{lpr}
  \end{table}

In the experiments in Section \ref{srn} the
  parameters are set to $R=8$ (a stencil
of 9 nodes to perform a least squares extrapolation), $r=4$ (degree of
the least squares interpolating polynomial), $r_0=2$ (substencils of
size 3 where smoothness indicators are computed),
$s_1=r_0=2$ (requirement to match the whole scheme
accuracy at the
boundary provided that the first and second
  derivatives do not vanish simultaneously), $s_2=1$ and $m=2$ (in
order to mimic the exponent choice at
the WENO-JS weights computation for spatial biased reconstructions,
although in the boundary case it does not affect the order of the
extrapolation when a discontinuity crosses the stencil since the
alternative choice is just the constant extrapolation).

\section{Numerical experiments}\label{srn}

\subsection{One-dimensional experiments}
In this section we present some one-dimensional
  numerical experiments where both the accuracy of the extrapolation method 
  for smooth solutions and its behavior in  presence of discontinuities
   will be tested and analyzed.

\subsubsection{Linear advection,  $\mathcal{C}^{\infty}$ solution.}

In this experiment we illustrate the performance of the proposed extrapolation method for a smooth solution.
The problem statement for this test is the same as in \cite{TanShu}.
We consider the following problem for the linear advection equation
\begin{equation}\label{eq:advcinf1}
\left\{
\begin{array}{l}
\displaystyle{u_t+u_x=0, \quad  (x,t)\in (-1,1) \times \mathbb{R}^+},\\
\displaystyle{u(x,0)= 0.25+0.5\sin(\pi x), \quad x\in (-1,1)}, \\
\displaystyle{u(-1, t)=0.25-0.5\sin(\pi(1+t)), \quad t\geq0}.
\end{array}\right.
\end{equation}

Outflow numerical boundary conditions are set at $x=1$. 
The unique solution of this problem
is $$u(x,t)=0.25+0.5\sin(\pi(x-t)).$$

We perform tests at uniform meshes with resolutions given by $n=20 \cdot 2^j$ points,
$j=1,\dots,5$. 
As we use a fifth order accurate
spatial scheme extrapolation is performed at 3 ghost cells at each
side of the boundary
On the other hand, the time discretization is integrated with a 
third order method, hence  we select a time
step given by  $\Delta
t=\left(\frac{2}{n}\right)^{\frac{5}{3}}$, with corresponding Courant
numbers $\Delta t / h_x=(2/n)^{2/3}\leq 1/20^{2/3}$ in order to attain 
fifth order accuracy in both integrators. Additionally, 
as the left boundary conditions are time dependent,
a specific approximation is needed in each of the 3 stages of the RK3-TVD time integrator at the boundary in order not to loose accuracy, as indicated in \cite{Carpenter}. The  boundary values required in this case are the following:
\begin{itemize}
  \item First stage: $g(t_k).$
  \item Second stage: $g(t_k)+\Delta tg'(t_k).$
  \item Third stage: $g(t_k)+\frac{1}{2}\Delta tg'(t_k)+\frac{1}{4}\Delta
    t^2g''(t_k).$
\end{itemize}
We run the simulation until  $t=1$ for all the previously
specified resolutions and for different modalities of 
boundary extrapolation: IW, WLS-UW and WLS-GAW. 
We show the errors corresponding to $1-$ and $\infty-$norm, and the numerical order computed from them in Tables \ref{IW}--\ref{WLS-GAW}.
\begin{table}[htb]
  \centering
  \begin{tabular}{|c|c|c|c|c|}
    \hline
    $n$ & Error $\|\cdot\|_1$ & Order $\|\cdot\|_1$ & Error
    $\|\cdot\|_{\infty}$ & Order $\|\cdot\|_{\infty}$ \\
    \hline
    40 & 8.73E$-6$ & $-$ & 2.44E$-5$ & $-$  \\
    \hline
    80 & 2.70E$-7$ & 5.01 & 7.35E$-7$ & 5.05 \\
    \hline
    160 & 8.45E$-9$ & 5.00 & 2.31E$-8$ & 4.99 \\
    \hline
    320 & 2.64E$-10$ & 5.00 & 6.95E$-10$ & 5.06 \\
    \hline
    640 & 8.26E$-12$ & 5.00 & 2.13E$-11$ & 5.03 \\
    \hline
  \end{tabular}
  \caption{Error table for problem \eqref{eq:advcinf1}, $t=1$, IW.}
  \label{IW}
\end{table}
\begin{table}[htb]
  \centering
  \begin{tabular}{|c|c|c|c|c|}
    \hline
    $n$ & Error $\|\cdot\|_1$ & Order $\|\cdot\|_1$ & Error
    $\|\cdot\|_{\infty}$ & Order $\|\cdot\|_{\infty}$ \\
    \hline
    40 & 4.28E$-5$ & $-$ & 1.99E$-4$ & $-$  \\
    \hline
    80 & 5.32E$-7$ & 6.33 & 1.86E$-6$ & 6.74 \\
    \hline
    160 & 1.38E$-8$ & 5.26 & 4.65E$-8$ & 5.32 \\
    \hline
    320 & 4.16E$-10$ & 5.06 & 1.43E$-9$ & 5.02 \\
    \hline
    640 & 1.27E$-11$ & 5.03 & 4.49E$-11$ & 5.00 \\
    \hline
  \end{tabular}
    \caption{Error table for problem \eqref{eq:advcinf1}, $t=1$, WLS-UW.}
  \label{WLS-UW}
\end{table}
\begin{table}[htb]
  \centering
  \begin{tabular}{|c|c|c|c|c|}
    \hline
    $n$ & Error $\|\cdot\|_1$ & Order $\|\cdot\|_1$ & Error
    $\|\cdot\|_{\infty}$ & Order $\|\cdot\|_{\infty}$ \\
    \hline
    40 & 1.50E$-5$ & $-$ & 4.29E$-5$ & $-$  \\
    \hline
    80 & 4.56E$-7$ & 5.04 & 1.44E$-6$ & 4.90 \\
    \hline
    160 & 1.37E$-8$ & 5.06 & 4.57E$-8$ & 4.98 \\
    \hline
    320 & 4.16E$-10$ & 5.04 & 1.43E$-9$ & 5.00 \\
    \hline
    640 & 1.27E$-11$ & 5.03 & 4.49E$-11$ & 4.99 \\
    \hline
  \end{tabular}
      \caption{Error table for problem \eqref{eq:advcinf1}, $t=1$, WLS-GAW.}
  \label{WLS-GAW}
\end{table}

From the results, and comparing with those obtained in
  \cite{BaezaMuletZorio2015} it can be seen that IW behaves essentially as
Lagrange extrapolation even for low resolutions, as it happens with
thresholding extrapolation with not excessively restrictive parameter
for the detection of discontinuities.

On the other hand, the two remaining techniques involving least
squares extrapolation produce slightly less accurate results, but still fifth
order accurate. This is what should be expected since a
wider stencil
is used, involving a polynomial of the same degree than IW. Albeit
this fact, we will see in the next experiments that the WLS-X
techniques are more robust than the ones based on Lagrange
extrapolation for more demanding problems.

\subsubsection{Burgers equation.}

We now consider Burgers equation
$$u_t+\left(\frac{u^2}{2}\right)_x=0,\quad\Omega=(-1,1),$$
with initial condition $u(x,0)=0.25+0.5\sin(\pi x)$, left inflow
boundary condition given by $u(-1,t)=g(t)$ and outflow
condition at the right boundary, where $g(t)=w(-1,t)$,
with $w$ the exact solution of the problem using periodic boundary
conditions. Numerical errors for spatial
resolutions $n=40\cdot2^k,$ $0\leq k\leq 5$, corresponding to time $t=0.3$ are presented in Tables \ref{burgersiw}-\ref{burgersgaw}. The solution is smooth 
and all methods provide fifth order accuracy as expected. 

The conclusions for this test are the same than those drawn in the linear advection problem.

\begin{table}[htb]
  \centering
  \begin{tabular}{|c|c|c|c|c|}
    \hline
    $n$ & Error $\|\cdot\|_1$ & Order $\|\cdot\|_1$ & Error
    $\|\cdot\|_{\infty}$ & Order $\|\cdot\|_{\infty}$ \\
    \hline
    40 & 2.03E$-5$ & $-$ & 3.57E$-4$ & $-$  \\
    \hline
    80 & 6.56E$-7$ & 4.95 & 1.47E$-5$ & 4.60 \\
    \hline
    160 & 1.36E$-8$ & 5.59 & 2.81E$-7$ & 5.71 \\
    \hline
    320 & 2.82E$-10$ & 5.58 & 9.16E$-9$ & 4.94 \\
    \hline
    640 & 7.58E$-12$ & 5.22 & 2.76E$-10$ & 5.05 \\
    \hline
    1280 & 2.23E$-13$ & 5.09 & 8.30E$-12$ & 5.06 \\
    \hline
  \end{tabular}
  \caption{Error table for Burgers equation, $t=0.3$. IW.}
  \label{burgersiw}
\end{table}

\begin{table}[htb]
  \centering
  \begin{tabular}{|c|c|c|c|c|}
    \hline
    $n$ & Error $\|\cdot\|_1$ & Order $\|\cdot\|_1$ & Error
    $\|\cdot\|_{\infty}$ & Order $\|\cdot\|_{\infty}$ \\
    \hline
    40 & 1.88E$-3$ & $-$ & 3.65E$-2$ & $-$  \\
    \hline
    80 & 6.62E$-5$ & 4.82 & 3.94E$-3$ & 3.21 \\
    \hline
    160 & 1.72E$-7$ & 8.52 & 9.45E$-6$ & 8.70 \\
    \hline
    320 & 2.50E$-9$ & 6.10 & 9.45E$-6$ & 5.37 \\
    \hline
    640 & 4.49E$-11$ & 5.80 & 6.73E$-9$ & 5.08 \\
    \hline
    1280 & 9.61E$-13$ & 5.55 & 1.92E$-10$ & 5.13 \\
    \hline
  \end{tabular}
  \caption{Error table for Burgers equation, $t=0.3$. WLS-UW.}
  \label{burgersuw}
\end{table}

\begin{table}[htb]
  \centering
  \begin{tabular}{|c|c|c|c|c|}
    \hline
    $n$ & Error $\|\cdot\|_1$ & Order $\|\cdot\|_1$ & Error
    $\|\cdot\|_{\infty}$ & Order $\|\cdot\|_{\infty}$ \\
    \hline
    40 & 7.80E$-4$ & $-$ & 2.61E$-2$ & $-$  \\
    \hline
    80 & 2.82E$-6$ & 8.11 & 8.14E$-5$ & 8.32 \\
    \hline
    160 & 1.11E$-7$ & 4.66 & 6.31E$-6$ & 3.69 \\
    \hline
    320 & 2.40E$-9$ & 5.53 & 2.28E$-7$ & 4.79 \\
    \hline
    640 & 4.48E$-11$ & 5.74 & 6.73E$-9$ & 5.08 \\
    \hline
    1280 & 9.61E$-13$ & 5.54 & 1.92E$-10$ & 5.13 \\
    \hline
  \end{tabular}
  \caption{Error table for Burgers equation, $t=0.3$. WLS-GAW.}
  \label{burgersgaw}
\end{table}

\begin{figure}[htb]
\centering
\subfloat[$t=1.1$.]{\label{fb:1}\includegraphics[width=0.47\textwidth]{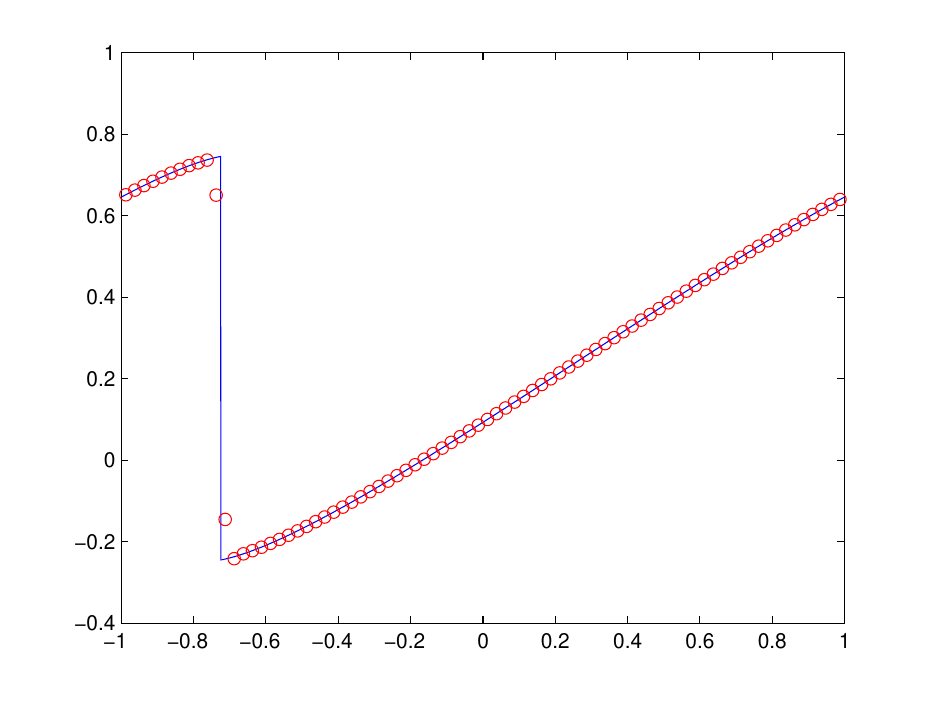}}
\subfloat[$t=12$.]{\label{fb:2}\includegraphics[width=0.47\textwidth]{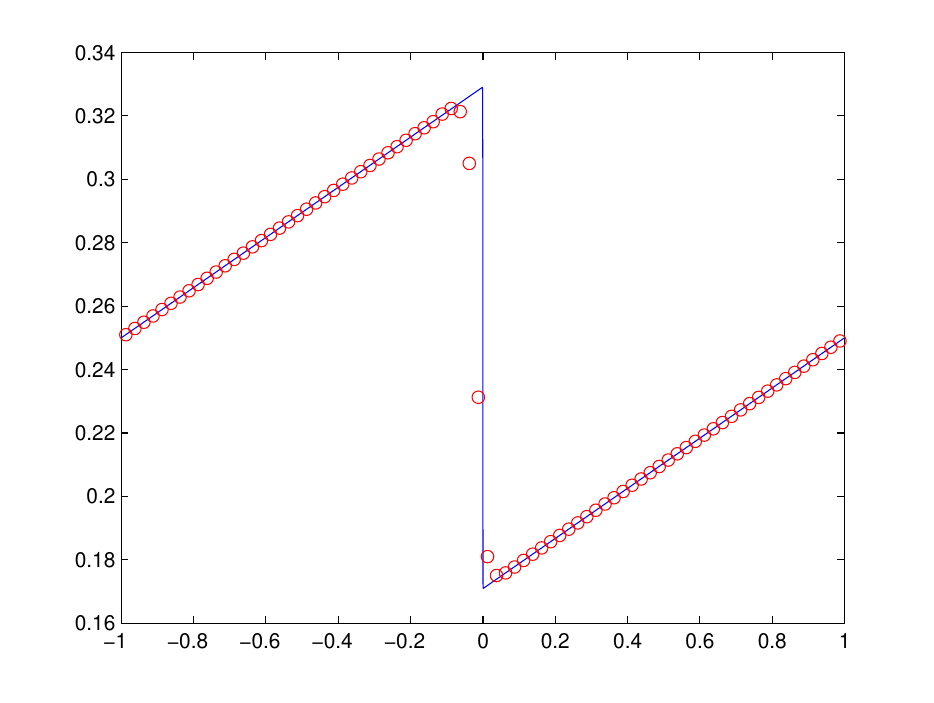}}
\caption{Shock in Burgers equation, $n=80$.}
\label{burgersshock}
\end{figure}

At further times a shock develops in the interior of the domain as 
shown in Figure \ref{fb:1}, corresponding to $t=1.1$; at
$t=8$ the shock interacts with the inflow boundary, being located at $x=0$
at $t=12$ as shown in Figure \ref{fb:2}.
Both plots in Figure \ref{burgersshock} illustrate that 
this discontinuity is properly captured by the scheme. 
Only the results for the WLS-GAW extrapolation method is included as the other methods provide similar results.

\subsubsection{Euler equations.}

The 1D Euler equations for gas dynamics are:

\begin{equation}\label{eq:eulereq1}
\begin{aligned}
&u_t+f(u)_x=0,\quad u=u(x,t),\quad\Omega=(0,1),\\
&u=\left[\begin{array}{c}
    \rho \\
    \rho v \\
    E \\
  \end{array}\right],\hspace{0.3cm}f(u)=\left[\begin{array}{c}
    \rho v \\
    p+\rho v^2 \\
    v(E+p) \\
  \end{array}\right],
\end{aligned}
\end{equation}
where $\rho$ is the density, $v$ is the
velocity and  $E$ is the specific energy of the system. The variable
$p$ stands for the pressure and is given by the equation of state:
$$p=\left(\gamma-1\right)\left(E-\frac{1}{2}\rho v^2\right),$$
where $\gamma$ is the adiabatic constant, that will be taken as
  $1.4$.

The initial data is the following, corresponding to the 
interaction of two blast waves:
$$u(x,0)=\left\{\begin{array}{ll}
    u_L & 0<x<0.1,\\
    u_M & 0.1<x<0.9,\\
    u_R & 0.9<x<1,
    \end{array}\right.$$
where $\rho_L=\rho_M=\rho_R=1$, $v_L=v_M=v_R=0$,
$p_L=10^3,p_M=10^{-2},p_R=10^2$. Reflecting boundary conditions are set
at $x=0$ and $x=1$, simulating a solid wall at both sides. This
problem involves multiple reflections of shocks and rarefactions off
the walls and many interactions of waves inside the domain. We will
use the same node setup as in the previous tests.

\begin{figure}[htb]
\centering
\subfloat[$\Delta x=1/800$.]{\label{fe:1}\includegraphics[width=0.47\textwidth]{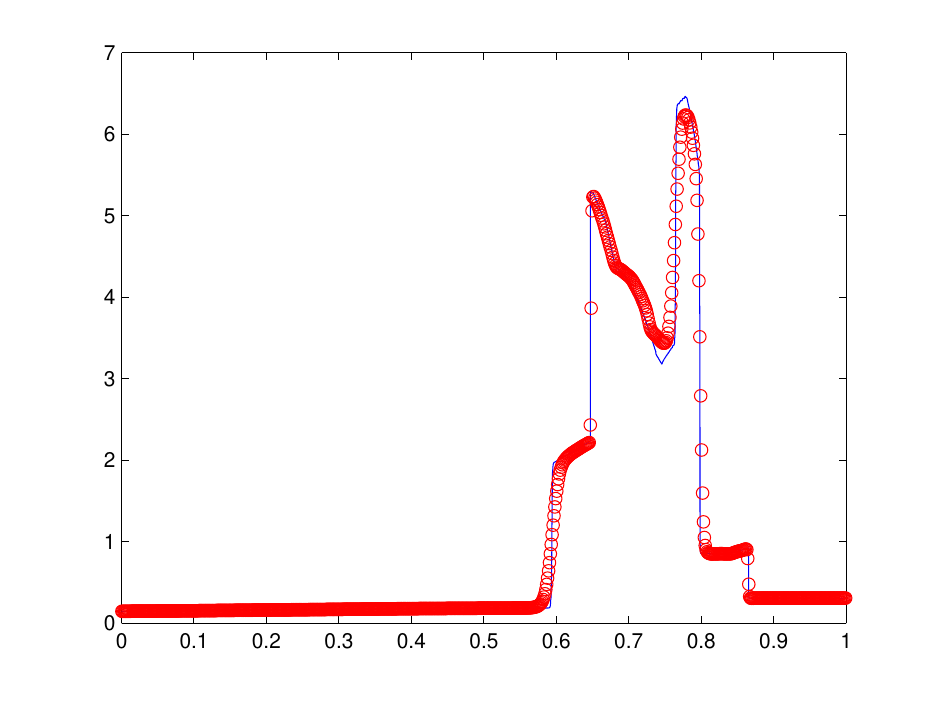}}
\subfloat[$\Delta x=1/1600$.]{\label{fe:2}\includegraphics[width=0.47\textwidth]{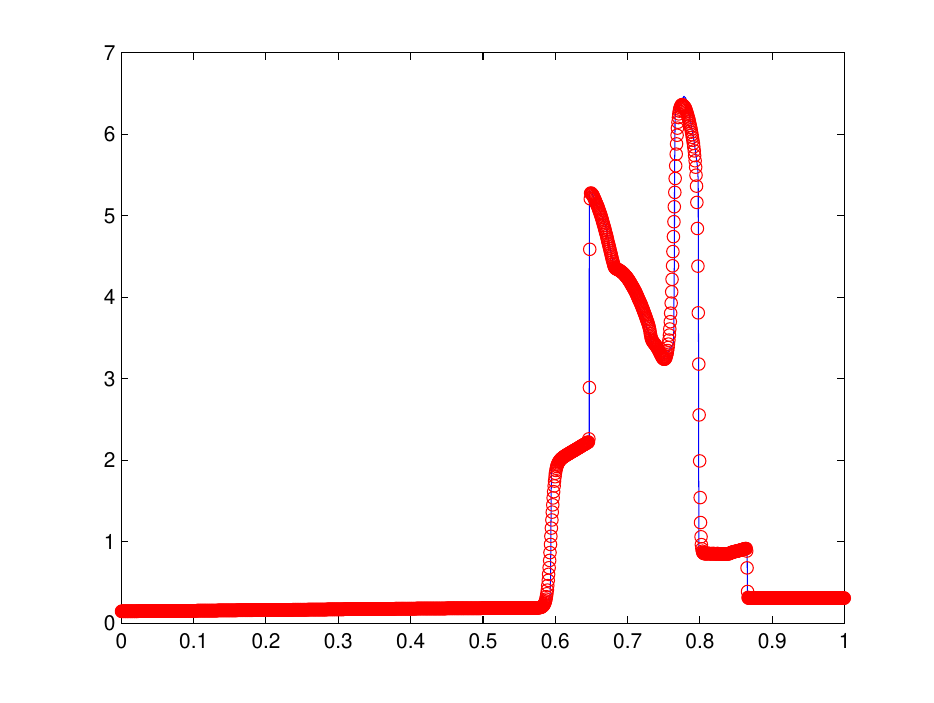}}
\caption{Density profile, $t=0.038$, WLS-GAW extrapolation.}
\label{eulershock}
\end{figure}

 Figure
\ref{eulershock} shows the density profile at $t=0.038$ at
two different resolutions, being the reference solution computed at
a resolution of $\Delta x=1/16000$. The
figure clearly shows that the results are satisfactory.
The plot corresponds to the WLS-GAW method being the results obtained with the other methods
almost identical

\subsection{Two-dimensional experiments}
\subsubsection{2D linear advection,  $\mathcal{C}^{\infty}$ solution.}
We consider the 2D linear advection equation
\begin{equation}\label{eq:adv2d1}
 u_t+u_x+u_y=0,
\end{equation}
with 
\begin{equation}\label{eq:adv2d2}
u_0=u(x,y,0)=0.25+0.5\sin(\pi(x+y))
\end{equation}
for different domains $\Omega$. We start with a square, $\Omega=(-1,1)$ where inflow
conditions
$g(t)=u(x,y,t)=0.25+0.5\sin(\pi(x+y-2t))$ are imposed at the left and
bottom boundary and outflow conditions at the rest. We compute the
solution with the same setup and techniques to achieve fifth order
accuracy as the above 1D tests for resolutions $n\times n$, for
$n=10\cdot 2^j,$ $1\leq j\leq 6$, and we obtain the results in Tables
\ref{S2DIW}-\ref{S2DWLS-GAW}
using different extrapolation techniques, running a simulation until
$t=1$.
\begin{table}[htb]
  \centering
  \begin{tabular}{|c|c|c|c|c|}
    \hline
    $n$ & Error $\|\cdot\|_1$ & Order $\|\cdot\|_1$ & Error
    $\|\cdot\|_{\infty}$ & Order $\|\cdot\|_{\infty}$ \\
    \hline
    20 & 5.36E$-4$ & $-$ & 1.56E$-3$ & $-$  \\
    \hline
    40 & 1.66E$-5$ & 5.01 & 5.40E$-5$ & 4.85  \\
    \hline
    80 & 5.24E$-7$ & 4.99 & 1.71E$-6$ & 4.98 \\
    \hline
    160 & 1.65E$-8$ & 4.99 & 5.33E$-8$ & 5.00 \\
    \hline
    320 & 5.15E$-10$ & 5.00 & 1.62E$-9$ & 5.04 \\
    \hline
    640 & 1.63E$-11$ & 4.98 & 5.04E$-11$ & 5.01 \\
    \hline
  \end{tabular}
  \caption{Error table for problem \eqref{eq:adv2d1} - \eqref{eq:adv2d2}, Square domain, IW.}
  \label{S2DIW}
\end{table}
\begin{table}[htb]
  \centering
  \begin{tabular}{|c|c|c|c|c|}
    \hline
    $n$ & Error $\|\cdot\|_1$ & Order $\|\cdot\|_1$ & Error
    $\|\cdot\|_{\infty}$ & Order $\|\cdot\|_{\infty}$ \\
    \hline
    20 & 6.14E$-3$ & $-$ & 9.07E$-2$ & $-$  \\
    \hline
    40 & 4.22E$-5$ & 7.19 & 2.98E$-4$ & 8.25  \\
    \hline
    80 & 6.90E$-7$ & 5.93 & 3.47E$-6$ & 6.42 \\
    \hline
    160 & 1.95E$-8$ & 5.15 & 9.13E$-8$ & 5.25 \\
    \hline
    320 & 5.94E$-10$ & 5.04 & 2.81E$-9$ & 5.02 \\
    \hline
    640 & 1.84E$-11$ & 5.01 & 8.90E$-11$ & 4.98 \\
    \hline
  \end{tabular}
    \caption{Error table for problem \eqref{eq:adv2d1} - \eqref{eq:adv2d2}, Square domain, WLS-UW.}
  \label{S2DWLS-UW}
\end{table}
\begin{table}[htb]
  \centering
  \begin{tabular}{|c|c|c|c|c|}
    \hline
    $n$ & Error $\|\cdot\|_1$ & Order $\|\cdot\|_1$ & Error
    $\|\cdot\|_{\infty}$ & Order $\|\cdot\|_{\infty}$ \\
    \hline
    20 & 8.22E$-4$ & $-$ & 2.07E$-3$ & $-$  \\
    \hline
    40 & 2.12E$-5$ & 5.28 & 8.10E$-5$ & 4.68  \\
    \hline
    80 & 6.39E$-7$ & 5.05 & 2.89E$-6$ & 4.81 \\
    \hline
    160 & 1.94E$-8$ & 5.04 & 9.03E$-8$ & 5.00 \\
    \hline
    320 & 5.94E$-10$ & 5.03 & 2.80E$-9$ & 5.01 \\
    \hline
    640 & 1.84E$-11$ & 5.01 & 8.91E$-11$ & 4.98 \\
    \hline
  \end{tabular}
    \caption{Error table for problem \eqref{eq:adv2d1} - \eqref{eq:adv2d2}, Square domain, WLS-GAW.}
  \label{S2DWLS-GAW}
\end{table}

\begin{figure}[htb]
  \centering
  \includegraphics[width=0.5\textwidth]{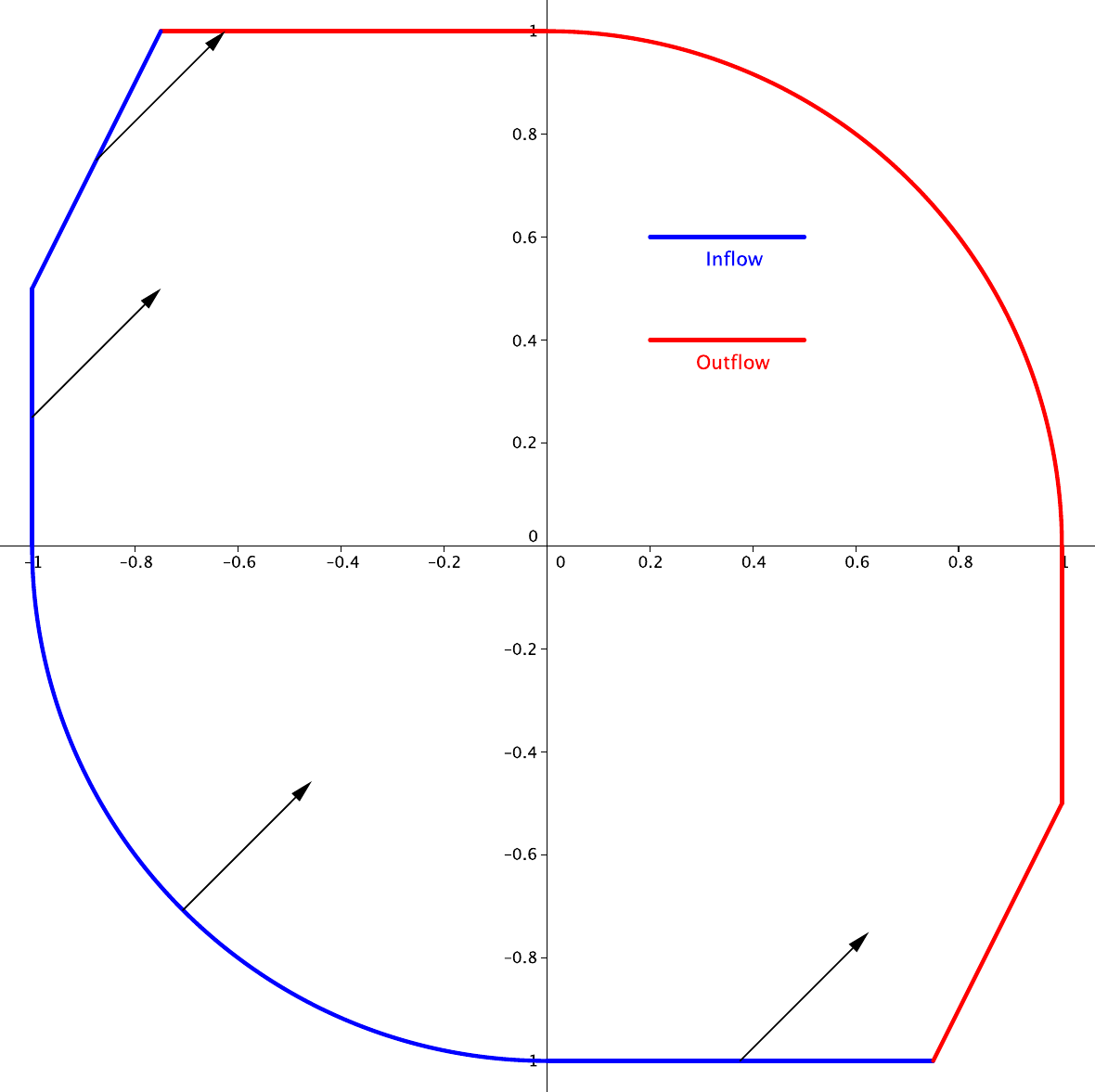}
  \caption{Complex domain for 2D experiment.}
  \label{fig:cd}
\end{figure}

We now change the domain and set $\Omega$ with respective boundary
conditions as indicated in Figure \ref{fig:cd}.

This time, the complexity of the domain makes Lagrange
extrapolation mildly unstable, whereas WLS is still stable and fifth
order accurate as
we can see numerically in Tables \ref{C2DWLS-UW}-\ref{C2DWLS-GAW} that
it is indeed achieved, running a simulation until $t=0.85$.

\begin{table}[htb]
  \centering
  \begin{tabular}{|c|c|c|c|c|}
    \hline
    $n$ & Error $\|\cdot\|_1$ & Order $\|\cdot\|_1$ & Error
    $\|\cdot\|_{\infty}$ & Order $\|\cdot\|_{\infty}$ \\
    \hline
    20 & 1.14E$-2$ & $-$ & 1.12E$-1$ & $-$  \\
    \hline
    40 & 4.82E$-3$ & 1.56 & 4.54E$-2$ & 1.40  \\
    \hline
    80 & 8.33E$-6$ & 9.18 & 2.25E$-4$ & 7.66 \\
    \hline
    160 & 7.53E$-8$ & 6.79 & 3.00E$-6$ & 6.23 \\
    \hline
    320 & 2.14E$-9$ & 5.14 & 1.01E$-7$ & 4.89 \\
    \hline
    640 & 6.59E$-11$ & 5.02 & 3.57E$-9$ & 4.82 \\
    \hline
  \end{tabular}
    \caption{Error table for problem \eqref{eq:adv2d1} - \eqref{eq:adv2d2}, Complex domain, WLS-UW.}
  \label{C2DWLS-UW}
\end{table}

\begin{table}[htb]
  \centering
  \begin{tabular}{|c|c|c|c|c|}
    \hline
    $n$ & Error $\|\cdot\|_1$ & Order $\|\cdot\|_1$ & Error
    $\|\cdot\|_{\infty}$ & Order $\|\cdot\|_{\infty}$ \\
    \hline
    20 & 3.78E$-3$ & $-$ & 3.06E$-2$ & $-$  \\
    \hline
    40 & 8.04E$-5$ & 5.56 & 1.41E$-3$ & 4.44  \\
    \hline
    80 & 1.81E$-6$ & 5.48 & 2.39E$-5$ & 5.88 \\
    \hline
    160 & 6.81E$-8$ & 4.73 & 2.42E$-6$ & 3.31 \\
    \hline
    320 & 2.13E$-9$ & 5.00 & 9.96E$-8$ & 4.60 \\
    \hline
    640 & 6.58E$-11$ & 5.01 & 3.57E$-9$ & 4.80 \\
    \hline
  \end{tabular}
      \caption{Error table for problem \eqref{eq:adv2d1} - \eqref{eq:adv2d2}, Complex domain, WLS-GAW.}
  \label{C2DWLS-GAW}
\end{table}

\subsubsection{Euler equations.}

The 2D Euler equations for inviscid gas dynamics are the following:
\begin{equation}\label{eq:eulereq}
\begin{aligned}
&u_t+f(u)_x+g(u)_y=0,\quad u=u(x,y,t),\\
&u=\left[\begin{array}{c}
    \rho \\
    \rho v^x \\
    \rho v^y \\
    E \\
  \end{array}\right],\hspace{0.3cm}f(u)=\left[\begin{array}{c}
    \rho v^x \\
    p+\rho (v^x)^2 \\
    \rho v^xv^y \\
    v^x(E+p) \\
  \end{array}\right],\hspace{0.3cm}g(u)=\left[\begin{array}{c}
    \rho v^y \\
    \rho v^xv^y \\
    p+\rho (v^y)^2 \\
    v^y(E+p) \\
  \end{array}\right].
\end{aligned}
\end{equation}
In these equations,  $\rho$ is the density, $(v^x, v^y)$  is the
velocity and  $E$ is the specific energy of the system. The variable
$p$  stands for the pressure and is given by the equation of state:
  $$p=(\gamma-1)\left(E-\frac{1}{2}\rho((v^x)^2+(v^y)^2)\right),$$
  where $\gamma$ is the adiabatic constant, that will be taken as
  $1.4$ in all the experiments.
\subsubsection{Double Mach Reflection}

This experiment models a vertical right-going Mach
10 shock colliding with an equilateral triangle. By symmetry, we run
the experiment in the upper half of the domain. The flow dynamics are
driven by the Euler equations \eqref{eq:eulereq} and the initial conditions
 are the following:
$$ u=(\rho,v^x,v^y,E)=\left\{\begin{array}{lcc}
 (8.0,8.25,0,563.5)&\textnormal{if}&x\leq \hat{x},\\
  (1.4,0,0,2.5)&\textnormal{if}&x>\hat{x}
\end{array}\right.
$$
where $\hat{x}$ is the point where the ramp starts.

\begin{figure}[htb]
  \centering
\begin{tabular}{cc}
  \includegraphics[width=0.47\textwidth]{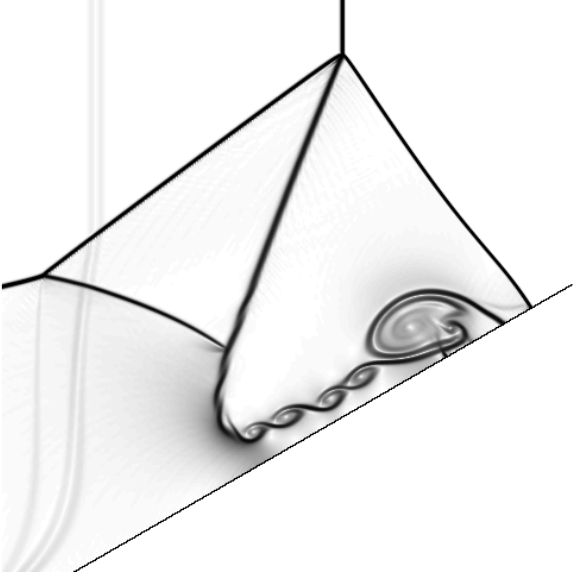}
  & \includegraphics[width=0.47\textwidth]{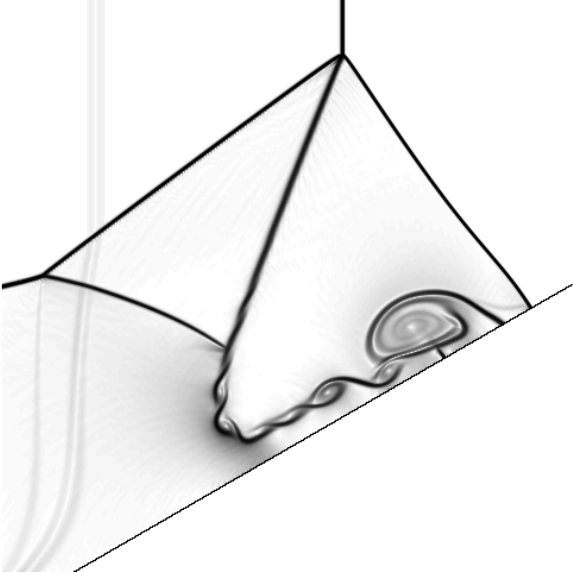} \\
  (a) First order extrapolation. & (b) WLS UW (WLS 0-UW). \\
  \includegraphics[width=0.47\textwidth]{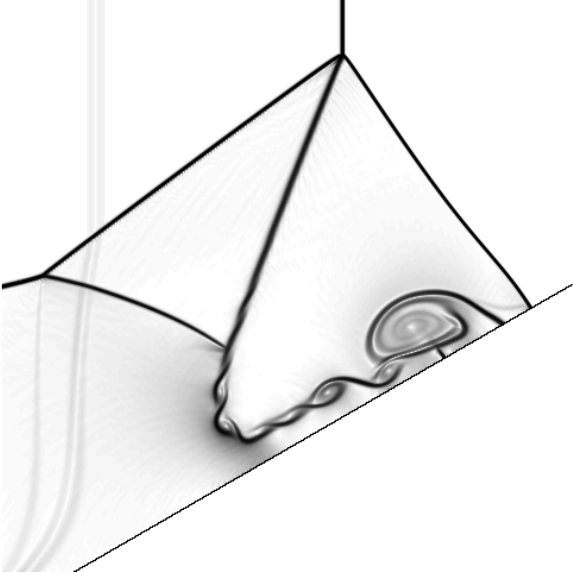}
  & \includegraphics[width=0.47\textwidth]{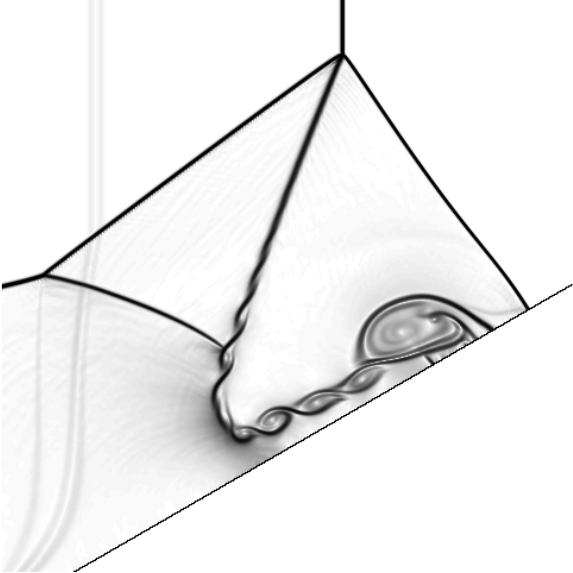} \\
  (c) WLS (-5)-UW & (d) WLS (-14)-UW. \\
  \includegraphics[width=0.47\textwidth]{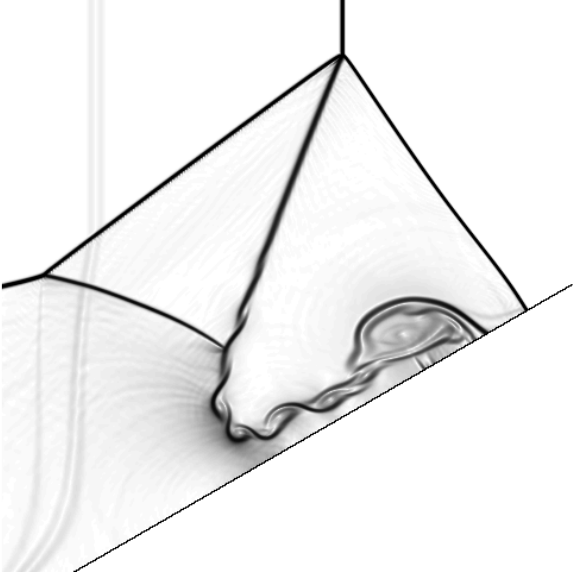}
  & \includegraphics[width=0.47\textwidth]{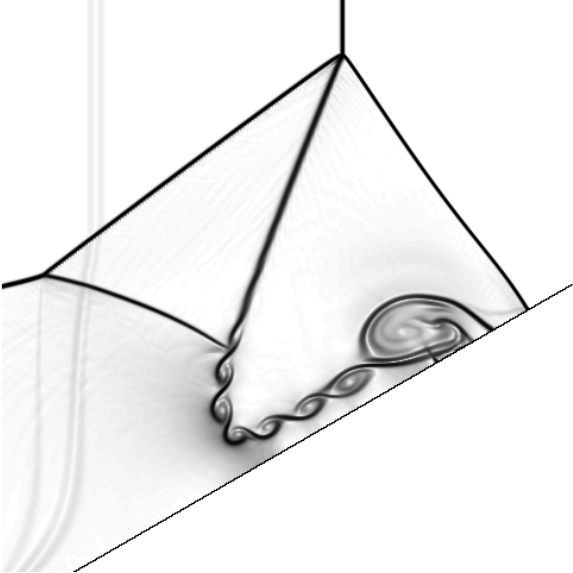} \\
  (e) WLS (-50)-UW & (f) WLS GAW. \\
\end{tabular}
\caption{Schlieren plots of the density for the Double Mach Reflection problem. $\displaystyle
  h_x=h_y=\frac{\sqrt{3}}{2}\frac{1}{512}$ (enlarged view)}
\label{DMR_512}
\end{figure}

\begin{figure}[htb]
  \centering
\begin{tabular}{cc}
  \includegraphics[width=0.47\textwidth]{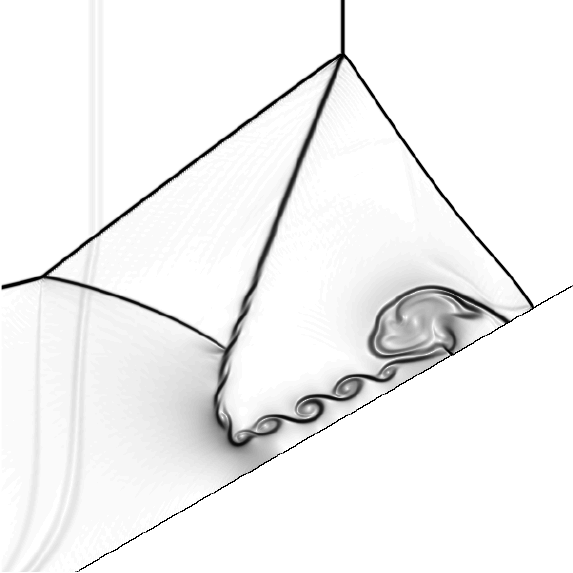}
  & \includegraphics[width=0.47\textwidth]{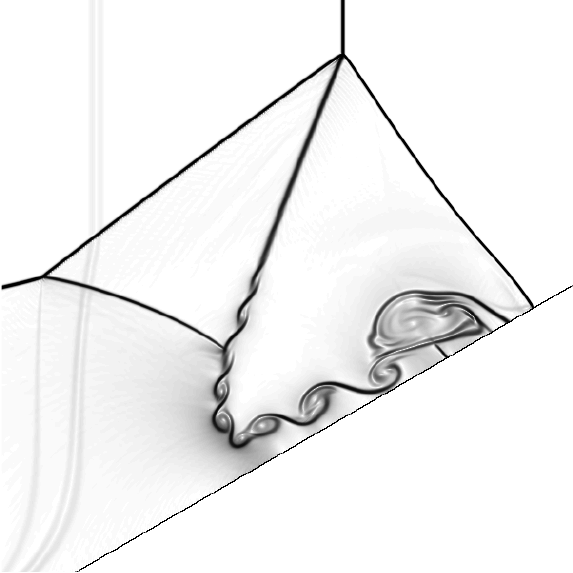} \\
  (a) First order extrapolation. & (b) WLS UW (WLS 0-UW). \\
  \includegraphics[width=0.47\textwidth]{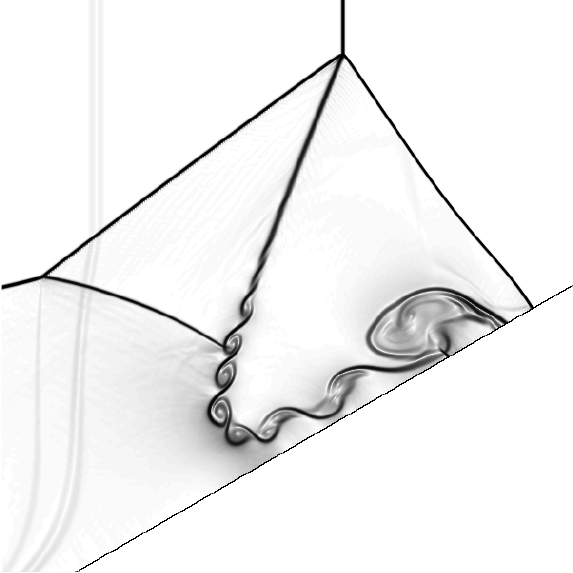}
  & \includegraphics[width=0.47\textwidth]{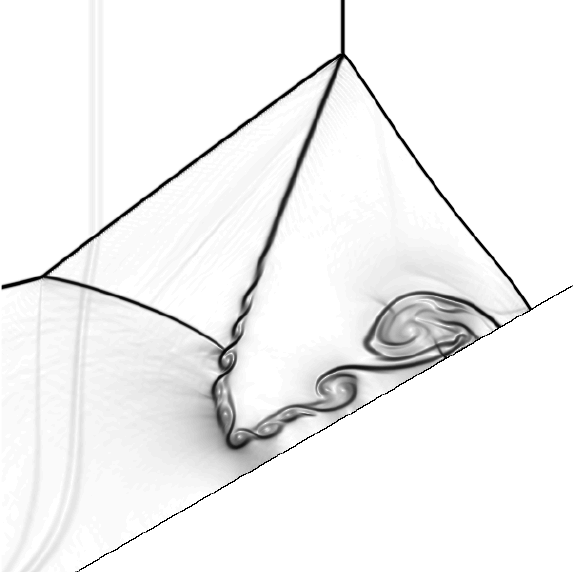} \\
  (c) WLS (-5)-UW & (d) WLS (-14)-UW. \\
  \includegraphics[width=0.47\textwidth]{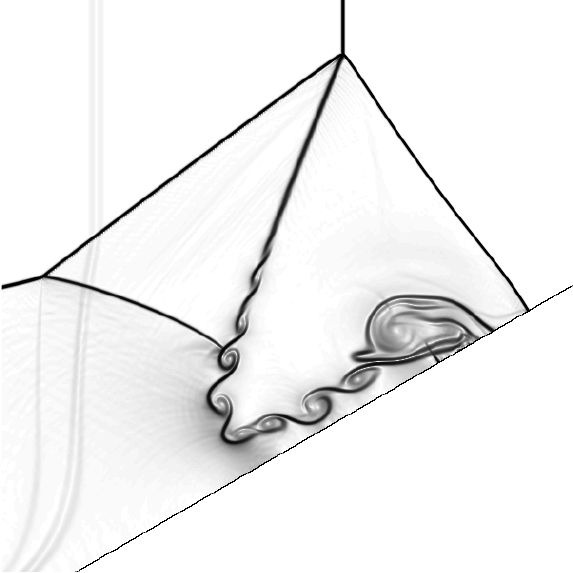}
  & \includegraphics[width=0.47\textwidth]{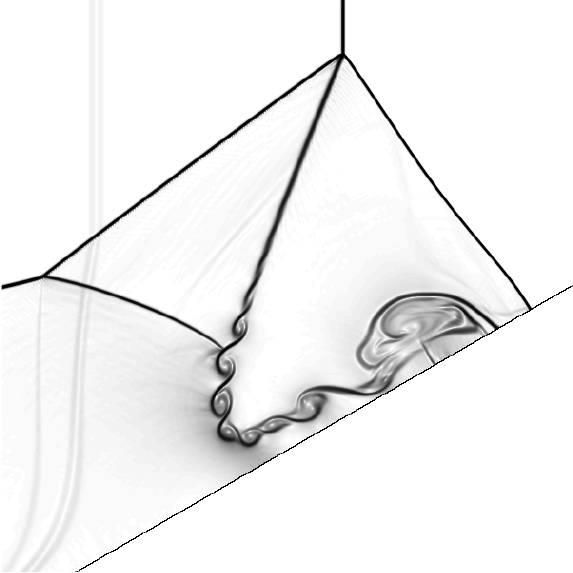} \\
  (e) WLS (-50)-UW & (f) WLS GAW. \\
\end{tabular}
\caption{Schlieren plots of the density for the Double Mach Reflection problem. $\displaystyle
  h_x=h_y=\frac{\sqrt{3}}{2}\frac{1}{640}$ (enlarged view)}
\label{DMR_640}
\end{figure}

We run the simulation for $t=0.2$. The experiment
consists in several simulations using WLS both for GAW 
and for $\lambda$-UW for various values of $\lambda$. 
 The results for two different resolutions are
presented as Schlieren plots of the density field in Figures \ref{DMR_512}-\ref{DMR_640}. The plots show the zone containing the
complex shock reflection structure.

It can be seen at the figures that the
performance of the new boundary extrapolation techniques are able to properly resolve the 
small-scale features present in this test for both methods, being the best results obtained with 
WLS-GAW and with WLS-$\lambda$-UW when taking the most negative values
for $\lambda$. Let us recall that the scheme will tend to
  achieve higher order extrapolations at the boundary as
  $\lambda\to-\infty$.

On the other hand we report in Figures \ref{DMR_512}(a) and \ref{DMR_640}(a) results 
obtained using a 
first order extrapolation procedure (based on a one-node stencil). It can be 
appreciated that the high order extrapolation techniques indeed entail an improvement with respect
to low order extrapolation, which illustrates the importance
  of keeping high order at the boundary even for scenarios where non-smooth structures appear.

\subsubsection{Interaction of a shock with multiple circular
  obstacles}
\begin{figure}[htb]
\centering
\begin{tabular}{cc}
  \raisebox{-0.35\textwidth}{
  \includegraphics[width=0.7\textwidth]{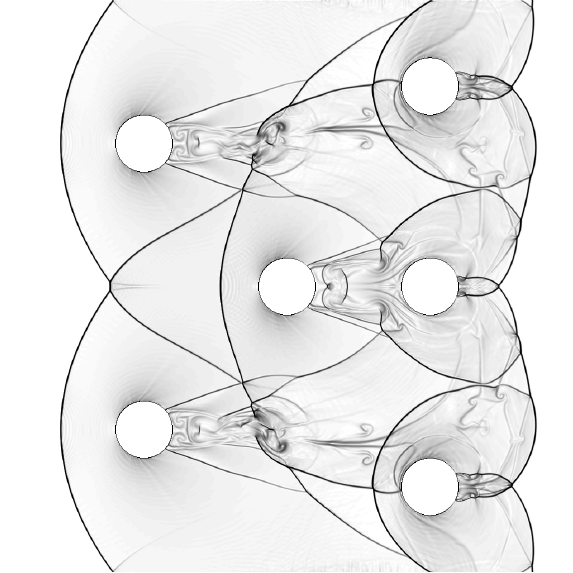}
  }
&
\begin{tabular}{@{}c@{}}
  \includegraphics[width=0.25\textwidth]{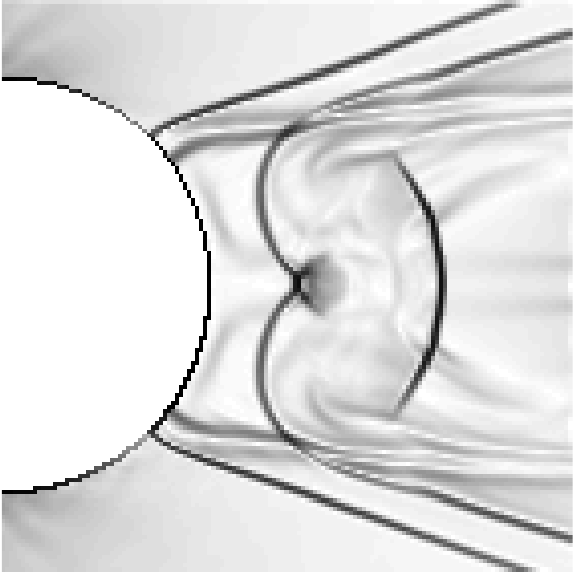}\\
  (a)\\[12.5ex]
  \includegraphics[width=0.25\textwidth]{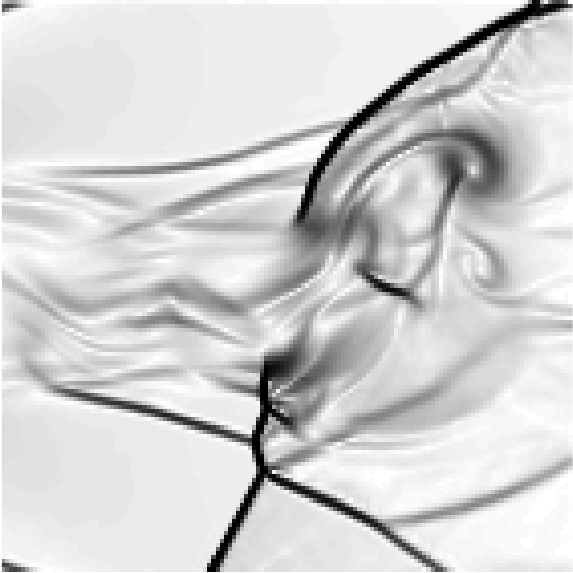}\\
  (b)
\end{tabular}
\end{tabular}
\caption{Circles. (a) Central circle (enlarged view). (b) Turbulences
  aside upper-left circle (enlarged view).}
\label{fig:circles}
\end{figure}

We now simulate a shock interacting with multiple solid circles. This
test can also be found in \cite{Boiron}, where penalization techniques 
were used to simulate problems including obstacles.
In this case, we run
the simulation until $t=0.5$ and a mesh size of
$h_x=h_y=\frac{1}{512}$ on the whole domain.  In this case, the
technique used for the numerical boundary conditions is WLS-GAW.
As in the previous experiment, we present a Schlieren plot for the last
time step in Figure \ref{fig:circles}.
These results are consistent with those obtained in \cite{Boiron}.
The result obtained in this experiment is very similar to the one
obtained through the thresholding technique in
\cite{BaezaMuletZorio2015}. The main difference is that this time a
thresholding parameter is no longer needed, which is an advantage
because in the case of thresholding extrapolation the performance
of this problem was strongly dependent on the choice of the parameter,
with moderate values leading to simulation failure and
requiring very high values in order to get satisfactory results.
\subsubsection{Steady-state supersonic flow around a triangle}
We next simulate the flow field over a solid triangle with height
$h=0.5$ and half angle $\theta=20$ deg moving at
supersonic speeds. The initial conditions are
$$u=(\rho,v^x,v^y,p)=(1, \sqrt{\gamma}M_1, 0, 1)$$
and the computation is stopped when a steady state is obtained from
the position of the shock waves. In order to halve the computational
cost, we perform the simulation in the upper half of the domain,
imposing appropriate reflecting boundary conditions at the
symmetry axis. We solve this problem using the WLS-GAW technique at
the boundary, as done in the previous experiment.

\begin{figure}[htb]
\centering
\begin{tabular}{cc}
\raisebox{-0.35\textwidth}{
  \includegraphics[width=0.7\textwidth]{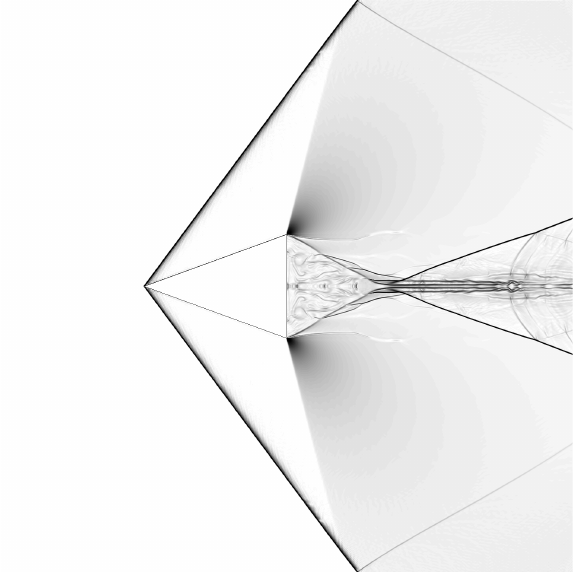}
  }
&
\begin{tabular}{@{}c@{}}
  \includegraphics[width=0.25\textwidth]{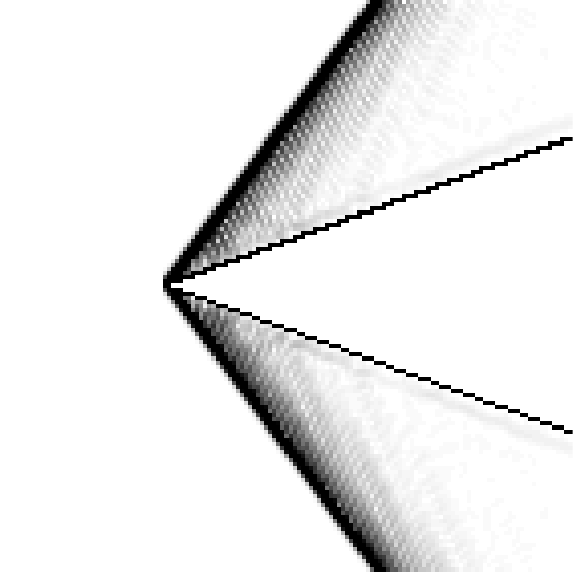}\\
  (a)\\[12.5ex]
  \includegraphics[width=0.25\textwidth]{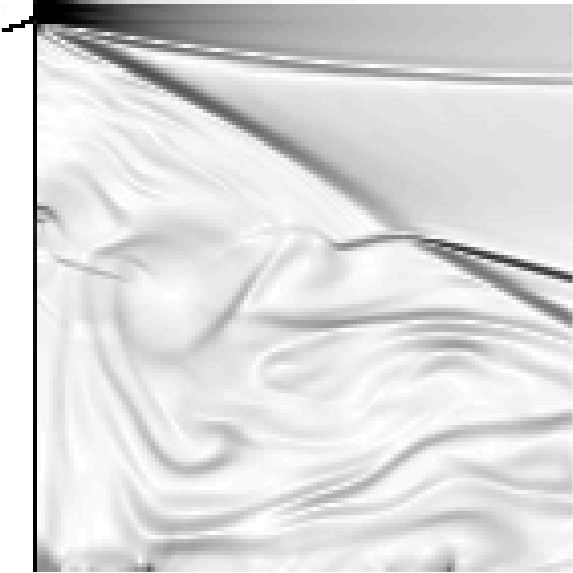}\\
  (b)
\end{tabular}
\end{tabular}
\caption{Triangle. (a) Left wedge (enlarged view). (b) Turbulences
  aside upper right corner of the triangle (enlarged view).}
\label{fig:triangle}
\end{figure}

Figure \ref{fig:triangle} shows
results that are consistent with those obtained in \cite{Boiron} and have sharper 
resolution that the ones reported in that paper for the same resolution.

\section{Conclusions and future work}\label{scn}
In this paper we have introduced some
  possible weighted
extrapolations akin to WENO reconstructions capable of keeping high
order accuracy for the global scheme,
which entails an improvement with respect
to the techniques based on thresholding parameters presented in
\cite{BaezaMuletZorio2015}. On the other hand, as stated in
\cite{TanShu}, straight Lagrange extrapolation may lead to a mildly
unstable scheme for multidimensional problems with some complex
domains.

We have seen that an appropriate and efficient option
  to achieve good
results both on smooth and non-smooth problems is to combine least
squares with an unique weight design to reduce to the constant
extrapolation (copying the value of the closest node) if there is a
discontinuity in the extrapolation stencil. The results obtained with
that technique are satisfactory and robust, with a weight design that, unlike
those defined in \cite{TanShu}, is dimensionless and scale
independent.

From the experiments, it can be concluded that the WLS-GAW technique
is better than WLS-UW without having to be scaled artificially
(``magnetize'' the weights to 1) in order to obtain a less diffusive
profile. We have seen that the results of WLS-GAW are better than
WLS-UW even for quite negative $\lambda$ values.

Since we now have a fully developed boundary extrapolation strategy,
our next purpose is to develop a parallelized AMR code \cite{BaezaMulet2006}, exploiting the
main benefits of using this global scheme in the above
terms. The extension of all the mentioned techniques
  to 3D is also under consideration.

\end{document}